\documentclass[11pt, reqno]{amsart}

\usepackage{amssymb}

\evensidemargin\oddsidemargin

\begin{document}
\baselineskip=18pt
\setcounter{page}{1}

\newcommand{\eqnsection}{
\renewcommand{\theequation}{\thesection.\arabic{equation}}
    \makeatletter
    \csname  @addtoreset\endcsname{equation}{section}
    \makeatother}
\eqnsection
   

\def\a{\alpha}
\def\B{{\bf B}}
\def\cA{{\mathcal{A}}} 
\def\cD{{\mathcal{D}}} 
\def\cG{{\mathcal{G}}} 
\def\cH{{\mathcal{H}}} 
\def\cL{{\mathcal{L}}} 
\def\cI{{\mathcal{I}}}
\def\CC{{\mathbb{C}}} 
\def\Dap{{\rm D}_{0+}^\a} 
\def\Dm{{\rm D}_{-}^\a} 
\def\Dp{{\rm D}_{+}^\a} 
\def\Ea{E_\a}
\def\esp{{\mathbb{E}}} 
\def\F{{\bf F}}
\def\Farl{{\F}_{\a,\lbd,\rho}}
\def\bF{{\bar F}}
\def\bG{{\bar G}}
\def\G{{\bf G}}
\def\Ga{{\Gamma}} 
\def\Gal{{\G_{\a,\lbd}}} 
\def\GG{{\bf \Gamma}}
\def\ii{{\rm i}} 
\def\Iap{{\rm I}_{0+}^\a} 
\def\Im{{\rm I}_{-}^\a} 
\def\Ip{{\rm I}_{+}^\a} 
\def\L{{\bf L}}
\def\lbd{\lambda}
\def\lacc{\left\{}
\def\lcr{\left[}
\def\lpa{\left(}
\def\lva{\left|}
\def\M{{\bf M}}
\def\NN{{\mathbb{N}}} 
\def\pb{{\mathbb{P}}}
\def\QQ{{\mathbb{Q}}} 
\def\R{{\bf R}}
\def\rl{{\mathbb{R}}}
\def\racc{\right\}}
\def\rpa{\right)}
\def\rcr{\right]}
\def\rva{\right|}
\def\sga{\sigma^{(\a)}}
\def\T{{\bf T}}
\def\Un{{\bf 1}}
\def\X{{\bf X}}
\def\Y{{\bf Y}}
\def\W{{\bf W}}
\def\Z{{\bf Z}}
\def\Warl{{\W}_{\a,\lbd,\rho}}
\def\Za{{\Z_\a}}

\def\elaw{\stackrel{d}{=}}
\def\claw{\stackrel{d}{\longrightarrow}}
\def\elaw{\stackrel{d}{=}}
\def\qed{\hfill$\square$}


\newtheorem{prop}{Proposition}[section]
\newtheorem{coro}[prop]{Corollary}
\newtheorem{lem}[prop]{Lemma}
\newtheorem{ex}[prop]{Example}
\newtheorem{exs}[prop]{Examples}
\newtheorem{rem}[prop]{Remark}
\newtheorem{theo}[prop]{Theorem}
\newtheorem{defi}[prop]{Definition}

\newcommand{\noi}{\noindent}
\newcommand{\dis}{\displaystyle }

\title{Fractional extreme distributions}

\author[L.~Boudabsa]{Lotfi Boudabsa}

\address{Institut de Math\'ematiques, Ecole Polytechnique F\'ed\'erale de Lausanne, CH-1015 Lausanne, Switzerland. {\em Email}: {\tt lotfi.boudabsa@epfl.ch}}

\author[T.~Simon]{Thomas Simon}

\address{Laboratoire Paul Painlev\'e, UMR 8524, Universit\'e de Lille,  Cit\'e Scientifique, F-59655 Villeneuve d'Ascq Cedex, France. {\em Email}: {\tt thomas.simon@univ-lille.fr}}

\author[P.~Vallois]{Pierre Vallois}

\address{Institut de Math\'ematiques Elie Cartan, INRIA-BIGS, Universit\'e de Lorraine, BP 239, F-54506 Vand\oe{}uvre-l\`es-Nancy Cedex, France. {\em Email}: {\tt pierre.vallois@univ-lorraine.fr}}

\keywords{Double Gamma function; Extreme distribution; Fractional differential equation; Kilbas-Saigo function; Le Roy function; Mittag-Leffler function; Peacock; Stable subordinator}

\subjclass[2010]{26A33; 33E12; 45E10; 60E05; 60E15; 60G52}

\begin{abstract} We consider three classes of linear differential equations on distribution functions, with a fractional order $\a\in [0,1].$ The integer case $\a =1$ corresponds to the three classical extreme families. In general, we show that there is a unique distribution function solving these equations, whose underlying random variable is expressed in terms of an exponential random variable and an integral transform of an independent $\a-$stable subordinator. From the analytical viewpoint, this law is in one-to-one correspondence with a Kilbas-Saigo function for the Weibull and Fr\'echet cases, and with a Le Roy function for the Gumbel case. By the stochastic representation, we can derive several analytical properties for the latter special functions, extending known features of the classical Mittag-Leffler function, and dealing with monotonicity, complete monotonicity, infinite divisibility, asymptotic behaviour at infinity, uniform hyperbolic bounds.
 
\end{abstract}

\maketitle

\section{Introduction and statement of the main results}

\label{intro}

The classical Fisher-Tippett-Gnedenko theorem states that the limit distributions arising from $a_n (\max(X_1, \ldots, X_n) - b_n)$ with $a_n > 0, b_n \in\rl$ and $(X_1, \ldots, X_n)$ an i.i.d real sample, can be classified up to positive affine transformation into three families:
$$\lacc\begin{array}{ll}
\W_\rho = \L^{1/\rho} & \quad \mbox{(Weibull distribution)}\\
\F_\rho = \L^{-1/\rho} & \quad \mbox{(Fr\'echet distribution)}\\
\G = \log \L & \quad \mbox{(Gumbel distribution)}
\end{array}\right.$$
where $\L$ is the unit exponential random variable, $\rho$ is a positive parameter and, with an abuse of notation which we will make throughout the paper, we have identified a random variable with its law. From the distribution function viewpoint, the three above extreme laws can also be obtained as the unique solution to a certain ordinary differential equation. More precisely, if $F(x)$ stands for a distribution function on $\rl$ and $\bF (x) = 1- F(x)$ denotes its associated survival function, the following equations
$$\lacc\begin{array}{lll}
F'(x) =  \rho x^{\rho -1} \bF(x), & x >0, & F(0) = 0\\
\bF'(x) = -\rho x^{-\rho -1} F(x), & x >0, & F(0) = 0\\
F'(x) = e^x \bF(x), & x\in\rl & 
\end{array}\right.$$
have each a unique solution which is respectively given by the Weibull distribution $\W_\rho,$ the Fr\'echet distribution $\F_\rho$ and the Gumbel distribution $\G.$ 

Notice that those three equations involve a logarithmic derivative and that they are solved via the exponential function. In this paper, we will consider some extensions of these equations in the context of fractional calculus. Throughout, we shall refer to the Appendix for all definitions and notations on the fractional integrals and derivatives that we will consider. In fractional calculus, a fundamental role is played by the classical Mittag-Leffler function
$$E_\a(z) \; =\; \sum_{n\ge 0} \frac{z^n}{\Ga(1+n\a)}, \qquad \a > 0,\, z\in\CC,$$which can be viewed as a generalization of the exponential function. We refer to Chapter 3 in \cite{GKMR} for a modern account on this function, and also to Chapter 2 therein for an interesting historical overview.

Let us first discuss an example. It is well-known from the general results of Barrett \cite{B1954} - see also Lemma 3.24 and the inversion formula (E.1.10) in \cite{GKMR} - that for every $\a,\lbd > 0$ the function $x\mapsto \Ea(-\lbd x^\a)$ solves on $\rl^+$ the following fractional differential equation
\begin{equation}
\label{Barr}
\Dap (1 -f)(x)\; =\;  \lbd f(x),
\end{equation}
where $\Dap$ is the progressive Liouville fractional derivative on the half-axis. Besides, it follows from the works of Pillai \cite{P1990} that for every $\a\in(0,1]$ the function $x\mapsto \Ea(-x^\a)$ is the survival function of a distribution on $\rl^+.$ More precisely, one has
\begin{equation}
\label{MLU}
\Ea(-x^\a)\; =\; \esp\lcr e^{-x^\a \Z_\a^{-\a}}\rcr\; =\; \pb\lcr \Za \times \L^{\frac{1}{\a}} > x\rcr
\end{equation}
where $\Za$ has a standard positive $\a-$stable distribution with the normalization $\esp[e^{-x\Za}] = e^{-x^\a}$ and, here and throughout, the product is assumed to be independent. This shows that the distribution function of the random variable $\Za \times \L^{\frac{1}{\a}}$ solves the fractional differential equation
\begin{equation}
\label{Barre}
\Dap F (x)\; =\;  \bF(x)
\end{equation}
on $(0,\infty),$ with the initial condition $F(0) = 0.$ 

The above fact can be used to display another, dual example involving the regressive Liouville fractional derivative on the half-axis $\Dm$ and the generalized Mittag-Leffler function
$$E_{\a,\beta}(z) \; =\; \sum_{n\ge 0} \frac{z^n}{\Ga(\beta +n\a)}, \qquad \a,\beta > 0,\; z\in\CC.$$
On the one hand, it follows indeed from the formulas (4.10.13) and (E.2.6) in \cite{GKMR} that for every $\a,\lbd > 0$ the function $x\mapsto \Ga(\a)\, E_{\a,\a}(-\lbd x^{-\a})$ solves on $(0,\infty)$ the fractional differential equation 
$$\Dm (1 -f)(x)\; =\;  \lbd x^{-2\a} f(x).$$
On the other hand, the above Pillai result shows that for every $\a\in(0,1],$ the function
$$\Ga(\a)\, E_{\a,\a}(-x^{-\a})\; =\; \Ga(1+\a) \,\Ea'(-x^{-\a})\; =\; \pb\lcr  (\Z_\a^{-1})^{(\a)}\times\L^{-\frac{1}{\a}}\le x\rcr$$
is a distribution function on $\rl^+,$ where the second equality follows from an elementary transformation of (\ref{MLU}) involving the size-bias $(\Z_\a^{-1})^{(\a)}$ of order $\a$ of the inverse positive $\a-$stable random variable (recall that for $t\in\rl$ and a positive random variable $X$ such that $\esp [X^t] < \infty$ the size-bias of order $t$ of $X$ is the random variable $X^{(t)}$ whose law is defined by
$$\esp [f(X^{(t)})]\; =\; \frac{\esp [ X^t f(X)]}{\esp[X^t]}$$
for all $f :\rl^+\to\rl$ bounded continuous, and that $\esp[\Z_\a^{-\a}] = 1/\Ga(1+\a)$). All of this shows that the distribution function of the random variable $(\Z_\a^{-1})^{(\a)} \times \L^{-1/\a}$ solves the fractional differential equation
\begin{equation}
\label{Barres}
\Dm \bF (x)\; =\;  x^{-2\a} F(x)
\end{equation}
on $(0,\infty),$ with the initial condition $F(0) = 0.$\\
 
In this paper, we wish to study more general fractional equations than (\ref{Barre}) and (\ref{Barres}), which are natural extensions of the above differential equations characterizing the classical extreme distributions.  Our findings involve the $\a-$stable subordinator $\{\sga_t\, t\ge 0\},$ which is the real L\'evy process starting from zero such that $\sga_1 \elaw \Za.$ For every $\a\in [0,1],$ its Laplace transform is given by
$$\esp[e^{-\lbd \sga_t}]\; =\; e^{-t\lbd^\a},\qquad \lbd, t\ge 0.$$  
Observe that $\sga$ is a pure drift for $\a = 1$ that is $\sigma^{(1)}_t =t,$ and a pure killing at an exponential time $\L$ for $\a =0$ that is $\sigma^{(0)}_t = \infty$ on $\{t \ge \L\}$ and $\sigma^{(0)}_t = 0$ on $\{t < \L\}.$ Our first main result gives a fractional extension of the Weibull distribution.

\begin{theo} 
\label{Weib}
For every $\lbd, \rho > 0$ and $\a\in [0,1],$ there exists a unique distribution function solving the fractional differential equation
\begin{equation}
\label{Weib+}
\Dap F (x)\; =\; \lbd\, x^{\rho -\a} \bF (x)
\end{equation}
on $(0,\infty),$ with the initial condition $F(0) = 0.$ The corresponding random variable is
$$\Warl \; \elaw\; \W_\rho \; \times\; \lpa \lbd \int_0^\infty \lpa (1- \sga_t)_+^{}\rpa^{\rho -\a} \, dt\rpa^{-\frac{1}{\rho}}.$$
\end{theo}

In the above statement, we have used the standard notation $x_+ = \max (x,0)$ for $x\in\rl.$ Observe that the integral on the right-hand side is finite a.s. for every $\rho > 0$ and $\a\in [0,1]:$ this is clear for $\a =1$ since $\rho > 0,$ and when $\a < 1$ this follows from the fact that $\sga$ is a non-decreasing c\`adl\`ag process which crosses the level 1 a.s. by a jump - see e.g. Theorem III.4 in \cite{B1996}. The above result shows that the fractional index $\a\in [0,1]$ of the derivative $\Dap$ gives rise to a non-trivial multiplicative perturbation of the Weibull random variable $\W_\rho$ given by the power of a certain Riemannian integral of the stable subordinator, whereas the parameter $\lbd$ is simply a scaling constant with $\Warl \elaw \lbd^{-1/\rho}\, \W_{\a, 1,\rho}.$ One has also the identities\begin{equation}
\label{WeibIDs}
\W_{1, \rho,\rho}\; \elaw\; \W_\rho\qquad\quad\mbox{and}\qquad\quad \W_{0, 1,\rho}\; \elaw\; \frac{\W_\rho}{\W_\rho}
\end{equation}
where, here and throughout, the quotient is assumed to be independent. The random variable on the right has a Pareto distribution of type III - see \cite{ARY} for a study of the latter distribution, and the mapping in law $\a\mapsto \W_{\a,\rho^\a\!, \rho}$ can be viewed as a parametrized arc connecting this Pareto III distribution and the Weibull distribution, the parameter being the index of the underlying stable subordinator. \\

Our second main result gives a fractional extension of the Fr\'echet distribution.

\begin{theo} 
\label{Fresh}
For every $\lbd, \rho > 0$ and $\a\in [0,1],$ there exists a unique distribution function solving the fractional differential equation 
\begin{equation}
\label{Fresh+}
\Dm \bF(x)\; =\; \lbd\, x^{-\rho -\a} F(x)
\end{equation}
on $(0,\infty),$ with the initial condition $F(0) = 0.$ The corresponding random variable is
$$\Farl \; \elaw\; \F_\rho \; \times\; \lpa \lbd \int_0^\infty \lpa 1+ \sga_t\rpa^{-\rho -\a} \, dt\rpa^{\frac{1}{\rho}}.$$
\end{theo}

In the above statement the integral on the right-hand side is finite a.s. by the law of the iterated logarithm at infinity for $\sga$ - see e.g. Theorem III.11 in \cite{B1996}. As above, the index $\a$ of the derivative $\Dm$ produces a multiplicative perturbation of the Fr\'echet random variable $\F_\rho$ via a Riemannian integral of $\sga$, whereas the parameter $\lbd$ is a scaling constant with $\Farl \elaw \lbd^{1/\rho}\, \F_{\a, 1,\rho}.$ One has also the identities
\begin{equation}
\label{FreshIDs}
\F_{1, \rho,\rho}\; \elaw\; \F_\rho\qquad\quad\mbox{and}\qquad\quad \F_{0, 1,\rho}\; \elaw\; \frac{\F_\rho}{\F_\rho}\; \elaw\; \frac{\W_\rho}{\W_\rho}\cdot
\end{equation}
It is interesting to observe from the two above theorems that the two mappings in law 
$$\a\;\mapsto\; \L^{\pm\frac{1}{\rho}} \; \times\;\lpa \rho^\a  \int_0^\infty \lpa 1\mp \sga_t\rpa^{\pm\rho -\a}_+ \, dt\rpa^{\mp\frac{1}{\rho}},$$
with a sign switch at $\a = 0$ corresponding to the trivial equation $F = x^{\rho} \bF,$ produce a parametrized arc connecting the classical extreme random variables $\W_\rho$ and $\F_\rho.$ The traditional role of the $\a-$stable subordinator is to define a fractional Laplacian via the underlying subordinated semi-group whose marginals are the symmetric $\beta-$stable distributions with $\beta = 2\a$, the densities of the latter being up to some multiplicative constant the solutions to the fractional Cauchy problem
$$\frac{\partial f}{\partial t} \; =\; \frac{\partial^{\beta} f}{\partial x^{\beta}} $$
on $(0,\infty)\times\rl^+$ - see \cite{GM1998} and the references therein for more on this standard subject. The above results show that the $\a-$stable subordinator is also involved, by means of its Riemannian integrals, in the solution to some fractional differential equations naturally associated to the Weibull or the Fr\'echet distribution. \\

The classical Gumbel distribution is the limit in law of either $\rho(\W_\rho -1)\claw\G$ or $\rho(1-\F_\rho)\claw\G$ as $\rho\to\infty.$ In order to define a fractional Gumbel distribution, it is natural from the above to introduce the random variable
$$\G_\a \; =\; \log\lpa\int_0^\infty e^{-\sga_t}\, dt\rpa$$
for every $\a\in [0,1],$ the a.s. convergence of the integral being a well-known consequence of  
$$\esp\lcr\int_0^\infty e^{-\sga_t}\, dt\rcr\; =\;\int_0^\infty \esp\lcr e^{-\sga_t}\rcr\, dt\; =\;\int_0^\infty e^{-t}\, dt\; =\;1. $$ 
We then have the following result involving the progressive Liouville fractional derivative on the line $\Dp$, which can be guessed at the formal limit $\rho\to\infty$ after a logarithmic change of variable in Theorem \ref{Weib} or Theorem \ref{Fresh}, and whose derivation is actually rigorous.

\begin{theo} 
\label{Gumb}
For every $\lbd > 0$ and $\a\in [0,1],$ there exists a unique distribution function solving the fractional differential equation
\begin{equation}
\label{Gumb+}
\Dp F (x)\; =\; \lbd^\a e^{\lbd x}\, \bF(x)
\end{equation}
on $\rl.$ The corresponding random variable is
$$\Gal \; \elaw\;\lbd^{-1}\lpa \G \, -\, \G_{\a}\rpa.$$
\end{theo}

Unlike the fractional Weibull or Fr\'echet distributions, the perturbation on the standard Gumbel distribution induced by the parameter $\a$ of the progressive derivative $\Dp$ is linear and not multiplicative. Again, the parameter $\lbd$ is a scaling constant with $\Gal \elaw \lbd^{-1} \G_{\a,1}.$ One has also the identities
$$\G_{1, 1}\; \elaw\; \G\qquad\quad\mbox{and}\qquad\quad \G_{0, 1}\; \elaw\; \G\, -\,\G.$$
The random variable on the right has a standard logistic distribution - see \cite{St} for an account on the latter distribution, and the mapping in law $\a\mapsto \G_{\a,1}$ can be viewed as a parametrized arc connecting the logistic distribution and the Gumbel distribution. \\

The proof of the three above theorems is divided in two parts. In Section 2 we prove the uniqueness of the solutions to a more general class of equations. These results have an independent interest, because the uniqueness problem is not always addressed in the literature on fractional calculus. It is classical in analysis to show that the solution of a differential equation must be the fixed point of an integral equation, and we use the same method, in the framework of fractional calculus. In Section 3 we show the existence of a probability law solving the above equations, and we establish the explicit multiplicative, respectively additive, factorizations. This is done via a one-to-one correspondence with a Kilbas-Saigo function, respectively a Le Roy function, which leads to a family of positive random variables characterized by their entire moments and previously studied in \cite{LS}, in a more general context. \\

The Kilbas-Saigo functions are three-parameter generalizations of the classical Mittag-Leffler functions $\Ea$ and $E_{\a,\beta}$ defined by the convergent series representation
$$E_{\a,m,l} (z)\; =\; \sum_{n\ge 0} \lpa\prod_{k=1}^n \frac{\Ga(1+ \a((k-1)m +l))}{\Ga(1 + \a((k-1)m + l +1))}\rpa z^n, \qquad z\in\CC,$$
for $\a,m >0$ and $l > -1/\a,$ with the convention made here and throughout that an empty product always equals 1. Note that $E_\a = E_{\a,1,0}$ and that $\Ga(\beta) E_{\a,\beta} = E_{\a,1,\frac{\beta-1}{\a}}.$ We refer to Chapter 5.2 in \cite{GKMR} for an account, including an extension to complex values of the parameter $l.$ The Le Roy functions are simple generalizations of the exponential function defined by
$$\cL_\a(z)\; =\; \sum_{n\ge 0} \frac{z^n}{(n!)^\a}, \qquad z\in\CC,$$
for $\a > 0.$ Introduced in \cite{L1899} in the context of analytic continuation, these functions are much less studied than Mittag-Leffler's. See however the recent paper \cite{GMR} and the references therein. In this paper, we can hinge upon the fractional extreme distributions to deduce some analytical features of these two interesting classes of special functions, in analogy to some known properties of the classical Mittag-Leffler functions. More precisely, we characterize their complete monotonicity on the negative half-line, we prove certain monotonicity properties with respect to the parameters, we derive their exact asymptotic behaviour at $-\infty$, and we establish uniform and optimal hyperbolic bounds. In particular, we prove the complete monotonicity of the function $x\mapsto E_{\a,m,m-1}(-x)$ for every $\a\in (0,1]$ and $m >0,$ solving an open question stated in \cite{DMV}. In a less complete way, we also study the infinite divisibility of the fractional extreme distributions. All these analytical results are to be found in Section 4. 

\section{Some uniqueness results on fractional hazard rates}

In this section we prove the uniqueness of distribution functions solving fractional equations of the type (\ref{Weib+}), (\ref{Fresh+}) or (\ref{Gumb+}), where the power function is replaced by a more general hazard rate.  We repeat that all definitions and notations on the fractional operators that we will consider here can be found in Appendix A.1.

\subsection{The Weibull case} We consider the equation
\begin{equation}
\label{pb1b+}
\Dap F\; =\; h \bF, \qquad F(0) =0,
\end{equation}
where $F$ is a distribution function and $h : (0,\infty)\to\rl^+$ is measurable and locally bounded. In the case $\a = 1,$ there exists a solution to (\ref{pb1b+}) if and only if
$$\int_0^\infty h(t)\, dt \; =\; \infty,$$ 
with a unique solution given by 
$$\bF(x)\; =\; \exp\Big\{-\int_0^xh(t)dt\Big\}\cdot$$
Recall that the function $h$ is called either the reliability function or the hazard rate of the underlying positive random variable. In the case $\a =0,$ there exists a solution to (\ref{pb1b+}) if and only if $h$ is non-decreasing, $h(0) = 0$ and $h(x)\to\infty$ as $x\to\infty,$ with a unique solution given by
$$\bF(x)\; =\; \frac{1}{1+h(x)}\cdot$$
In order to state our result in the fractional case $\a\in (0,1),$ we introduce the following linear operator
$$A^{\a,h}_{0+} : f\;\mapsto\; {\rm I}^{\alpha}_{0+}(hf)$$
which is well-defined on measurable functions from $(0,\infty)$ to $\rl^+.$ 

\begin{theo}\label{t1+}
Assume that there exists $\rho > 0$ such that $h(x) = O(x^{\rho-\a})$ as $x\to 0.$ Then, if it exists, the distribution function satisfying \eqref{pb1b+} is uniquely defined  by the convergent series 
\begin{equation}
\label{rpf6+}
\bF (x) \; =\; \sum_{n\ge 0} (-1)^n (A^{\a,h}_{0+})^n\Un\, (x), \qquad x > 0.
\end{equation}
\end{theo}

\proof

We begin by transforming \eqref{pb1b+} into an integral equation. Since $F$ is a distribution function on $\rl^+,$ there exists a probability measure $\mu$ on $\rl^+$ such that
\begin{equation}
\label{ItFr}
F(x)\; =\; \int_0^x \mu(dt)\; =\; {\rm I}^{\alpha}_{0+}(F_{\a,\mu})(x),\qquad x > 0,
\end{equation}
where we have set
$$F_{\a,\mu}(x) \; =\; \frac{1}{\Ga(1-\a)}\int_0^x (x-t)^{-\a} \,\mu(dt), \qquad x >0,$$
and the second equality in (\ref{ItFr}) is a direct consequence of Fubini's theorem and of the evaluation of the Beta integral of the first kind. Moreover, it is easy to see that the function $F_{\a,\mu},$ which might take infinite values, is nevertheless locally integrable at zero since $\mu$ is a probability. Hence, applying ${\rm I}^\a_{0+}$ on both sides of (\ref{pb1b+}), we can use the inversion formula (\ref{Fr0+}) and get
$$F\; =\; {\rm I}^{\alpha}_{0+}( h\bF)\; =\; A^{\a,h}_{0+}(\bF)$$
on $(0,\infty).$ This leads to the fixed point equation $\bF = \Un - A^{\a,h}_{0+}(\bF)$ and, by the linearity of $A^{\a,h}_{0+}$, to
$$\bF (x) \; =\; \sum_{k = 0}^{n-1} (-1)^k (A^{\a,h}_{0+})^k\Un\, (x)\; +\; (-1)^{n} (A^{\a,h}_{0+})^{n}(\bF)\, (x), \qquad x > 0,$$
for every $n\ge 1.$ Fixing now $x >0,$ the assumption made on $h$ implies that there exists a constant $c >0$ such that $h(t) \le c t^{\rho -\a}$ for every $t\in [0,x].$ Since moreover $\bF(t)\le 1$ for every $t\in[0,x],$ an immediate induction based on the Beta integral of the first kind implies
$$0\;\le\;(A^{\a,h}_{0+})^{n}(\bF)\, (x)\; \le\; \lpa\prod_{k=1}^n \frac{\Ga(1-\a +k\rho)}{\Ga(1+k\rho)}\rpa c^n x^{\rho n}\;\to\; 0\quad\mbox{as $n\to\infty,$}$$
where the convergence towards zero follows e.g. from (1.1.5) in \cite{AAR1999}. This completes the proof.

\endproof

\begin{rem}
\label{rt1+}
{\em (a) It can be proved that the series defined in \eqref{rpf6+} converges uniformly on every compact set and that it tends to 1 as $x\to 0.$ On the other hand, because of the alternate signs it is difficult to guess whether it is non-increasing and converges to zero as $x\to \infty.$ In the particular case when $h$ is a power function, Theorem \ref{Weib} gives a positive answer with the help of a Kilbas-Saigo function and the $\a-$stable subordinator. It would be interesting to know if there are other hazard rate functions $h$ such that the series in \eqref{rpf6+} is indeed a survival function on $(0,\infty).$

\medskip

(b)  In a different direction, sharing a certain analogy with the previous item,  the authors have introduced in \cite{TaVa2016, TaVa2017} generalized fractional distributions which are not conventional and classical distributions with fractional hazard rates. In \cite{TaVa2018}, the stochastic approximation of fractional probability distribution have been studied.

\medskip

(c) The above proof shows the more general fact that under the same assumption on $h$, for every $\a\in (0,1)$ there exists a unique bounded solution to
$$\Dap (f(0) -f)\; =\;  h f$$
on $(0,\infty),$ which is given by
$$f\; =\; f(0)\;\times\;\sum_{n\ge 0} (-1)^n (A^{\a,h}_{0+})^n\Un.$$
This can be viewed as an extension of (\ref{Barr}) which handles the case when $h$ is a positive constant.}

\end{rem}

\subsection{The Fr\'echet case} We consider the equation
\begin{equation}
\label{pb1b-}
\Dm \bF\; =\; h F, \qquad F(0) =0,
\end{equation}
given on distribution functions, where $h : (0,\infty)\to\rl^+$ is measurable and locally bounded. In the case $\a = 1,$ there exists a solution to (\ref{pb1b-}) if and only if $h$ is locally integrable on $(0, \infty]$ and such that 
$$\int_0^\infty h(t)\, dt \; =\; \infty,$$ 
with a unique solution given by
$$F(x)\; =\; \exp\Big\{-\int_x^\infty h(t)dt\Big\}\cdot$$
 In the case $\a =0,$ there exists a solution to (\ref{pb1b-}) if and only if $h$ is non-increasing, $h(0+) =\infty$ and $h(x)\to 0$ as $x\to \infty,$ with a unique solution given by
$$F(x)\; =\; \frac{1}{1+h(x)}\cdot$$
In order to state our result in the fractional case $\a\in (0,1),$ we introduce the following linear operator
$$A^{\a,h}_{-} : f\;\mapsto\; {\rm I}^{\alpha}_{-}(hf),$$
which is well-defined on measurable functions from $(0,\infty)$ to $\rl^+.$ 

\begin{theo}
\label{t1-}
Assume that there exists $\rho > 0$ such that $h(x) = O(x^{-\rho-\a})$ as $x\to \infty.$ Then, if it exists, the distribution function satisfying \eqref{pb1b-} is uniquely defined by the convergent series 
\begin{equation}
\label{rpf6-}
F (x) \; =\; \sum_{n\ge 0} (-1)^n (A^{\a,h}_{-})^n\Un\, (x), \qquad x > 0.
\end{equation}
\end{theo}

\proof

The proof is analogous to that of \eqref{t1-}, except that we deal with survival functions. Since $\bF$ is a survival function on $\rl^+,$ there exists a probability measure $\mu$ on $\rl^+$ such that
$$\bF(x)\; =\; \int_x^\infty \mu(dt)\; =\; {\rm I}^{\alpha}_{-}(\bF_{\a,\mu})(x),\qquad x > 0,$$
where we have set
$$\bF_{\a,\mu}(x) \; =\; \frac{1}{\Ga(1-\a)}\int_x^\infty (t-x)^{-\a} \,\mu(dt), \qquad x >0,$$
and used Fubini's theorem together with the evaluation of the Beta integral of the second kind. The function $\bF_{\a,\mu}$ is locally integrable at infinity since for every $y >0$ one has
$$\int_y^\infty \bF_{\a,\mu}(x) \, dx \; =\; \frac{1}{\Ga(2-\a)}\int_y^\infty (t-y)^{1-\a} \,\mu(dt) \; =\; \frac{1}{\Ga(2-\a)}\int_y^\infty (s-y)^{-\a} \bF(s) \,ds\; <\; \infty,$$
the finiteness of the third integral following from the fact that $\Dm\bF$ must be finite on $(0,\infty).$ Hence, we can apply the inversion formula (\ref{Fr0-}) and get
$\bF = {\rm I}^{\alpha}_{-}( h F) = A^{\a,h}_{-}(F)$ on $(0,\infty),$ which leads to $F = \Un - A^{\a,h}_{-}(F)$ and then to
$$F (x) \; =\; \sum_{k = 0}^{n-1} (-1)^k (A^{\a,h}_{-})^k\Un\, (x)\; +\; (-1)^{n} (A^{\a,h}_{-})^{n}(F)\, (x), \qquad x > 0,$$
for every $n\ge 1.$ Fixing now $x >0,$ the assumption made on $h$ implies that there exists a constant $c >0$ such that $h(t) \le c t^{-\rho -\a}$ for every $t > x.$ Since moreover $F(t)\le 1$ for every $t >x,$ an induction based on the Beta integral of the second kind implies
$$0\;\le\;(A^{\a,h}_{-})^{n}(F)\, (x)\; \le\; \lpa\prod_{k=1}^n \frac{\Ga(k\rho)}{\Ga(\a+k\rho)}\rpa c^n x^{-\rho n}\;\to\; 0\quad\mbox{as $n\to\infty,$}$$
which completes the proof.

\endproof

\begin{rem}
\label{rt1-}
{\em (a) It can be proved that the series defined in \eqref{rpf6-} converges uniformly on every compact set of $(0,\infty)$ and that it tends to 1 as $x\to \infty.$ As above, it is not easy to guess from the alternate signs whether the series is non-decreasing and converges to zero as $x\to 0.$ The case when $h$ is a power function gives a positive answer in Theorem \ref{Fresh}, with the help of another Kilbas-Saigo function. It would be interesting to know if there are other functions $h$ such that the series in \eqref{rpf6-} is indeed a distribution function on $(0,\infty).$

\medskip

(b) The above proof shows that under the same assumption on $h$, for every $\a\in (0,1)$ there exists a unique bounded function having a limit $\ell$ at infinity and solving
$$\Dm (\ell -f)\; =\;  h f$$
on $(0,\infty),$ which is given by
$$f\; =\; \ell \;\times\;\sum_{n\ge 0} (-1)^n (A^{\a,h}_{-})^n\,\Un.$$
Observe that this solution is zero if $\ell =0.$}
\end{rem}

\subsection{The Gumbel case} We consider the equation
\begin{equation}
\label{pb1b}
\Dp F\; =\; h \bF, \qquad F(x) > 0\;\; \mbox{on $\rl$}
\end{equation}
given on distribution functions, where $h : \rl\mapsto\rl^+$ is measurable and locally bounded. In the case $\a = 1,$ there exists a solution to (\ref{pb1b}) if and only if $h$ is locally integrable at $-\infty$ and such that 
$$\int_\rl h(t)\, dt \; =\; \infty,$$ 
with a unique solution given by
$$\bF(x)\; =\; \exp\Big\{-\int_{-\infty}^xh(t)dt\Big\}\cdot$$
In the case $\a =0,$ there exists a solution to (\ref{pb1b+}) if and only if $h$ is non-decreasing, $h(x) \to 0$ as $x\to -\infty$ and $h(x)\to\infty$ as $x\to\infty,$ with a unique solution given by
$$\bF(x)\; =\; \frac{1}{1+h(x)}\cdot$$
In order to state our result in the fractional case $\a\in (0,1),$ we introduce the following linear operator
$$A^{\a,h}_{+} : f\;\mapsto\; {\rm I}^{\alpha}_{+}(hf),$$
which is well-defined on measurable functions from $\rl$ to $\rl^+.$ The following result is a simple variation on Theorem \ref{t1-}.
 
\begin{theo}
\label{t1G}
Assume that there exists $\lbd > 0$ such that $h(-x) = O(e^{-\lbd x})$ as $x\to \infty.$ Then, if it exists, the distribution function satisfying \eqref{pb1b} is uniquely defined by the convergent series 
\begin{equation}
\label{rpf6}
\bF (x) \; =\; \sum_{n\ge 0} (-1)^n (A^{\a,h}_{+})^n\Un\, (x), \qquad x\in\rl.
\end{equation}
\end{theo}

\proof
Setting $G (x) = \bF(-x),$ we see that $G$ is a distribution function on $\rl$ with $\bG(x) = F(-x).$ Changing the variable transforms \eqref{pb1b} into
$$\Dm \bG\; =\; g G, \qquad G(x) > 0\;\; \mbox{on $\rl$},$$
with $g(x) = h(-x).$ An argument similar to the proof of Theorem \ref{t1-}, where the evaluation of the Beta integral of the second kind is replaced by that of the Gamma integral, gives then the unique solution
$$G (x) \; =\; \sum_{n\ge 0} (-1)^n (A^{\a,g}_{-})^n\Un\, (x), \qquad x \in\rl.$$
Changing the variable backwards, we obtain the required \eqref{rpf6}.

\endproof

\section{Proof of the main theorems}

In this section we show the existence of the real random variables associated to the fractional differential equations (\ref{Weib+}), (\ref{Fresh+}) and (\ref{Gumb+}), and we express them in terms of the unit exponential random variable and an integral transform of an independent $\a-$stable subordinator. The main ingredient in the proof is the following infinite independent product
$$\T(a,b,c)\; =\; \prod_{n\ge 0} \lpa\frac{a+nb +c}{a+nb}\rpa \B_{a+ nb, c}$$
where, here and throughout, $\B_{a,b}$ denotes a standard Beta random variable. We refer to Section 2.1 in \cite{LS} for more details on this product, including the fact that it is a.s. convergent for every $a,b,c > 0.$ We also mention from Proposition 2 in \cite{LS} that its Mellin transform is 
$$\esp[\T(a,b,c)^s]\;=\; \lpa \frac{\Ga(ab^{-1})}{\Ga ((a+c)b^{-1})}\rpa^s\times\; \frac{[a+c; b]_s}{[a; b]_s}$$
for every $s > -a,$ where $[z;\delta]_s$ stands for the generalized Pochhammer symbol which is defined in (\ref{Poch}) below. In general, we shall refer to the Appendix A.2 for all the properties of Barnes' double Gamma function and its associated Pochhammer symbol that we will need. 
 
\subsection{Proof of Theorem \ref{Weib}}
 
We first consider the case $\a\in (0,1).$ The uniqueness is a direct consequence of Theorem \ref{t1+} with $h(x) = \lbd x^{\rho -\alpha}.$ Moreover, we know by (\ref{rpf6+}) that a distribution solving (\ref{Weib+}), if it exists, has survival function
$$\bF (x) \; =\; \sum_{n\ge 0} (-1)^n (A^{\a,h}_{0+})^n\Un\, (x), \qquad x > 0.$$
Since
$$(\Iap t^\beta)(x) \; =\; \frac{\Ga(\beta + 1)}{\Ga(\a +\beta +1)} \, x^{\a +\beta +1}$$
for every $\beta > -1,$ an induction implies
$$\bF(x) \; = \; \sum_{n\ge 0}\left(\prod_{k=1}^{n}\frac{\Gamma(k\rho+1-\alpha)}{\Gamma(k\rho+1)}\right)(-\lambda  x^{\rho})^n\; = \; E_{\alpha,\frac{\rho}{\a},\frac{\rho}{\a}-1}(-\lambda x^\rho)$$
for every $x\ge 0.$ Observe that alternatively, the fact that $E_{\alpha,\frac{\rho}{\a},\frac{\rho}{\a}-1}(-\lambda x^\rho)$ is a solution to (\ref{Weib+}) follows from Theorem 5.27 and the inversion formula (E.1.10) in \cite{GKMR}.\\

It thus remains to prove that $x\mapsto E_{\alpha,\frac{\rho}{\a},\frac{\rho}{\a}-1}(-\lambda x^\rho)$ is a survival function on $\rl^+$ and to identify the underlying positive random variable. For every $z\in\CC,$ one has
$$E_{\alpha,\frac{\rho}{\a},\frac{\rho}{\a}-1}(z)\; =\; \sum_{n\ge 0} a_n(\a,\rho) \,\frac{z^n}{n!}$$
with 
$$a_n(\a,\rho)\;= \;\rho^{-n} \,\prod_{k=1}^n\frac{\Gamma(k\rho +1 -\alpha)}{\Gamma(k\rho)}\cdot$$
Let us now consider the positive random variable
$$\frac{\Ga(\rho +1-\a)}{\Ga(\rho +1)}\; \T(1,\rho^{-1}, (1-\a)\rho^{-1}).$$
By the aforementioned Proposition 2 in \cite{LS} and the concatenation formula (\ref{Conca1}), the positive entire moments of the latter random variable are given by
$$\delta^{n} \lpa \frac{[1+(1-\a)\delta;\delta]_n}{[1;\delta]_n}\rpa\; =\; a_n(\a,\rho)$$
where $\delta = \rho^{-1}.$ The Stirling formula (\ref{Stirling}) implies
$$a_n(\a,\rho)^{-\frac{1}{2n}}\; \sim\; \kappa n^{\frac{\a-1}{2}}\qquad 
\mbox{as $n\to \infty$}$$
for some positive constant $\kappa,$ so that Carleman's criterion
$$\sum_{n\ge 1} a_n(\a,\rho)^{-\frac{1}{2n}}\; =\; \infty$$
is fulfilled, and the law of the latter random variable is determined by its positive entire moments. Finally, it follows from Theorem (b) (i) in \cite{LS} with $q =\rho -\a$ that
$$\frac{\Ga(\rho +1-\a)}{\Ga(\rho +1)}\; \T(1,\rho^{-1}, (1-\a)\rho^{-1})\; \elaw\; \cA(\a, 0,\rho -\a)\; \elaw\; \int_0^\infty \lpa 1- \sga_t\rpa_+^{\rho -\a} \, dt.$$
All in all, we have shown that
\begin{equation}
\label{KS1}
E_{\alpha,\frac{\rho}{\a},\frac{\rho}{\a}-1}(z)\; =\; \esp\lcr\exp \lacc z\, \int_0^\infty \lpa 1- \sga_t\rpa_+^{\rho -\a} \, dt\racc\rcr, \qquad z\in\CC.
\end{equation}
This implies 
\begin{eqnarray*}
E_{\alpha,\frac{\rho}{\a},\frac{\rho}{\a}-1}(-\lambda x^\rho) & = & \esp\lcr\exp\lacc -\lbd x^{\rho}\, \int_0^\infty \lpa 1- \sga_t\rpa_+^{\rho -\a} \, dt\racc\rcr\\
& = & \pb\lcr\, \L\, >\, x^{\rho}\,\lpa \lbd\int_0^\infty \lpa 1- \sga_t\rpa_+^{\rho -\a} \, dt\rpa\rcr\; = \; \pb\lcr \Warl > x\rcr
\end{eqnarray*}
for every $x\ge 0,$ which completes the proof for $\a\in(0,1).$ The case $\a =1$ is that of the classical Weibull distribution and was already discussed in the introduction. Finally, the case $\a = 0$ amounts to solving $F = \lbd x^\rho (1-F),$ which yields $\bF(x) = 1/(1+ \lbd x^\rho)$ and
$$\W_{0,\lbd,\rho}\;\elaw\; \lpa \frac{\L}{\lbd\L}\rpa^{\frac{1}{\rho}}\;\elaw\; \W_\rho\,\times\,\lpa \lbd \int_0^\infty \lpa 1- \sigma^{(0)}_t\rpa_+^\rho \, dt\rpa^{-\frac{1}{\rho}}.$$
\qed

\begin{rem}
\label{rW1}
{\em (a) The main result of \cite{JSW2018} implies the identification
$$\lpa\frac{\Ga(\rho +1-\a)}{\Ga(\rho)}\; \T(1,\rho^{-1}, (1-\a)\rho^{-1})\rpa^{-\frac{1}{\rho}}\elaw\;\,\cG(\rho +1-\a, 1-\a)$$
where $\cG(m,a)$ is the generalized stable random variable with parameters $m > a > 0,$ whose density is up to normalization the unique positive solution to  
\begin{equation}
\label{JSW1}
 {\rm I}^a_{0+} f \; =\; x^m f 
\end{equation}
on $(0,\infty).$ This yields the further identity in law
\begin{equation}
\label{JSW2}
\Warl\; \elaw\; \lpa\rho\lbd^{-1}\rpa^{\frac{1}{\rho}}\W_\rho\,\times\, \cG(\rho +1-\a, 1-\a).
\end{equation}
In particular - see the introduction in \cite{JSW2018} for the third identity, one has
$$\W_{\a,1,\a}\; \elaw\; \a^{1/\a}\W_{\a,\a,\a}\; \elaw\;  (\a^{1/\a} \cG(1, 1-\a))\,\times\,\W_\a\; \elaw\; \Za\,\times\,\L^{1/\a},$$
in accordance with the Pillai distribution mentioned in the introduction since
$$\pb\lcr \W_{\a,1,\a} > x\rcr\; =\; E_{\alpha,1,0}(-x^\a)\; =\; E_{\a}(-x^\a), \qquad x\ge 0.$$ 
Notice that the integro-differential equation (\ref{JSW1}) shares some formal similarities with (\ref{Weib+}). Nevertheless it is essentially different because it deals with densities whereas (\ref{Weib+}) deals with distribution functions.

\medskip

(b) There exist unique solutions to fractional differential equations of the type (\ref{Weib+}) without the restriction to distribution functions. The main result of \cite{BDP2008} states that for every $\a\in (0,1],$ there is a unique solution to 
$$(\a+1)\Dap f\; = \; -x f$$
satisfying the boundary condition 
$$f(x)\; \sim\; \frac{x^{\a-1}}{\Ga(\a)\Ga(1/(\a+1))}$$ 
at zero, which is the density function of the running maximum of a spectrally positive $(\a+1)-$stable L\'evy process starting from zero.

\medskip

(c) With the above notation, one has
$$a_n(\a,\rho)\; =\; \frac{n!}{\Phi_{\a,\rho}(1)\cdots\Phi_{\a,\rho}(n)}$$
with
\begin{eqnarray*}
\Phi_{\a,\rho}(x) \; = \; \frac{\Ga(1+\rho x)}{\Ga(1-\a +\rho x)}& = & \frac{1}{\Ga(1-\a)}\; +\; \frac{\a}{\rho\Ga(1-\a)}\int_0^\infty (1- e^{-x u}) \frac{e^{-u\rho^{-1}}}{(1-e^{-u\rho^{-1}})^{\a +1}}\, du
\end{eqnarray*}
a Bernstein function. By Proposition 3.3 in \cite{CPY}, we deduce
\begin{equation}
\label{Txi}
\frac{\Ga(\rho +1-\a)}{\Ga(\rho+1)}\; \T(1,\rho^{-1}, (1-\a)\rho^{-1})\;\elaw\;\int_0^\infty e^{-\xi^{(\a,\rho)}_t}\, dt
\end{equation}
where $\lacc \xi^{(\a,\rho)}_t\! , \; t\ge 0\racc$ is the subordinator having Laplace exponent
$$\esp\lcr e^{-\lbd \xi^{(\a,\rho)}_t}\rcr\; = \; e^{-t \Phi_{\a,\rho}(\lbd)}.$$
This leads to
$$\Warl\; \elaw\; \W_\rho\;\times\; \lpa \lbd \int_0^\infty e^{-\xi^{(\a,\rho)}_t}\, dt \rpa^{-\frac{1}{\rho}}.$$
Let us notice that the identification (\ref{Txi}) can also be deduced from Corollary 5 in \cite{LS} is the case $q = \rho -\a > -\a$ and ${\hat \rho} =1,$ with the notation therein. Observe finally that this is consistent with the limiting case $\a = 1$ with $\xi^{(1,\rho)}_t = \rho\, t$ and $\a = 0$ with $\xi^{(0,\rho)}_t = \sigma^{(0)}_t.$

\medskip

(d) The above proof shows that the function
$$x\;\mapsto\;E_{\alpha,m,m-1}(-x)$$
is completely monotone (CM) for every $\a\in (0,1]$ and $m >0.$ This can be viewed as a generalization of the classic result by Pollard that $E_\a(-x) = E_{\alpha,1,0}(-x)$ is CM for every $\a\in (0,1]$ - see e.g. Proposition 3.23 in \cite{GKMR}. As already mentioned in the introduction, this also solves a conjecture stated in \cite{DMV} - see Section 4 and equations (10) and (11) therein. Notice that the formula (\ref{KS1}) implies the Bernstein representation
$$E_{\alpha,m,m-1}(-x)\; =\; \esp\lcr \exp -x \lacc \int_0^\infty \lpa 1- \sga_t\rpa_+^{\a (m-1)} \, dt\racc\rcr, \qquad x\ge 0.$$
For $m = 1,$ setting $T_\a = \inf\{ t > 0,\; \sga_t > 1\}\elaw \Z_\a^{-\a}$ we obtain
$$E_\a(-x)\; =\; E_{\alpha,1,0}(-x)\; =\; \esp [e^{-x T_\a}]\; =\; \esp [e^{-x \Z_\a^{-\a}}],$$
a well-known fact following from our discussion prior to the statement of Theorem \ref{Weib}. See \cite{TS2015} for other CM functions related to $E_\a.$ In Section \ref{FurthA} below, we will generalize this fact and show some further analytical properties of the Kilbas-Saigo function $E_{\alpha,m,m-1}.$}

\end{rem}

We end this section with a convergent series representation, in the non-explicit case $\a\in(0,1),$ for the density $f^\W_{\a,\lbd,\rho}$ of $\Warl.$ This is an immediate consequence of a term-by-term differentiation of the survival function
$$\pb[\Warl > x]\; =\; E_{\a,\frac{\rho}{\a},\frac{\rho}{\a} -1}(-\lbd x^\rho)$$
which was obtained during the proof of Theorem \ref{Weib}. 
 
\begin{coro}
\label{WeibD}
For every $\a\in (0,1),$ the density of $\Warl$ has the following convergent series representation on $(0,\infty)$
$$f^\W_{\a,\lbd,\rho}(x)\; =\; \lbd \, x^{\rho -1}\,\sum_{n\ge 0}\lpa \prod_{j=0}^n \frac{\Ga(j\rho + \rho + 1 -\a)}{\Ga(j\rho +\rho)}\rpa \frac{(-\lbd x^\rho)^n}{\rho^n\,n!}\cdot$$
\end{coro}

\subsection{Proof of Theorem \ref{Fresh}}

The proof is analogous to that of Theorem \ref{Weib}, except that we will deal with distribution functions instead of survival functions. We begin with the case $\a\in (0,1)$ and the uniqueness follows from Theorem \ref{t1-}. Besides, we know by (\ref{rpf6-}) that a distribution solving (\ref{Fresh+}), if it exists, must have distribution function
$$F (x) \; =\; \sum_{n\ge 0} (-1)^n (A^{\a,h}_{-})^n\Un\, (x), \qquad x > 0,$$
with $h(x) = \lbd x^{-\rho -\alpha}.$  Since
$$(\Im t^{-\beta})(x) \; =\; \frac{\Ga(\beta -\a)}{\Ga(\beta)} \, x^{\a -\beta}$$
for every $\beta > \a,$ an induction implies
$$F(x) \; = \; \sum_{n\ge 0}\left(\prod_{k=1}^{n}\frac{\Gamma(k\rho)}{\Gamma(k\rho+\a)}\right)(-\lambda  x^{-\rho})^n\; = \; E_{\alpha,\frac{\rho}{\a},\frac{\rho-1}{\a}}(-\lambda x^{-\rho})$$
for every $x > 0.$ Alternatively, the fact that $E_{\alpha,\frac{\rho}{\a},\frac{\rho-1}{\a}}(-\lambda x^{-\rho})$ is a solution to (\ref{Fresh+}) follows from Theorem 5.30 and the inversion formula (E.1.10) in \cite{GKMR}. It remains to prove that $x\mapsto E_{\alpha,\frac{\rho}{\a},\frac{\rho-1}{\a}}(-\lambda x^{-\rho})$ is a distribution function on $(0,\infty)$ and to identify the underlying positive random variable. For every $z\in\CC,$ one has
$$E_{\alpha,\frac{\rho}{\a},\frac{\rho-1}{\a}}(z)\; =\; \sum_{n\ge 0} b_n(\a,\rho) \,\frac{z^n}{n!}$$
with 
$$b_n(\a,\rho)\;= \;\rho^{-n} \,\prod_{k=1}^n\frac{\Gamma(k\rho +1)}{\Gamma(k\rho+\a)}\cdot$$
Reasoning exactly as above implies that $\{b_n(\a,\rho), \; n\ge 0\}$ is the determinate integer moment sequence of the positive random variable 
$$\frac{\Ga(\rho)}{\Ga(\rho +\a)}\; \T(1 +\a\rho^{-1},\rho^{-1}, (1-\a)\rho^{-1})\; \elaw\;\int_0^\infty \lpa 1+ \sga_t\rpa^{-\rho -\a} \, dt,$$
where the identity in law follows from Corollary 3 in \cite{LS}. We have hence shown that
\begin{equation}
\label{KS2}
E_{\alpha,\frac{\rho}{\a},\frac{\rho-1}{\a}}(z)\; =\; \esp\lcr\exp \lacc z\, \int_0^\infty \lpa 1+ \sga_t\rpa^{-\rho -\a} \, dt\racc\rcr, \qquad z\in\CC.
\end{equation}
This implies 
\begin{eqnarray*}
E_{\alpha,\frac{\rho}{\a},\frac{\rho-1}{\a}}(-\lambda x^{-\rho}) & = & \esp\lcr\exp\lacc -\lbd x^{-\rho}\, \int_0^\infty \lpa 1+ \sga_t\rpa^{-\rho -\a} \, dt\racc\rcr\\
& = & \pb\lcr\, \L\, \ge\, x^{-\rho}\,\lpa \lbd\int_0^\infty \lpa 1+ \sga_t\rpa^{-\rho -\a} \, dt\rpa\rcr\; = \; \pb\lcr \Farl \le x\rcr
\end{eqnarray*}
for every $x > 0,$ which completes the proof for $\a\in(0,1).$ As for Theorem \ref{Fresh} the remaining cases $\a =0$ and $\a =1$ are elementary and we leave the details to the reader.

\qed

\begin{rem}\label{rF1} {\em (a) Contrary to $\Warl,$ the factor
$$\int_0^\infty \lpa 1+ \sga_t\rpa^{-\rho -\a} \, dt$$
appearing in the decomposition of $\Farl$ cannot be expressed as a generalized stable law. On the other hand, this factor can also be viewed as the perpetuity of some subordinator: rewriting
$$b_n(\a,\rho)\; =\; \frac{n!}{\Psi_{\a,\rho}(1)\cdots\Psi_{\a,\rho}(n)}$$
with the Bernstein function
\begin{eqnarray*}
\Psi_{\a,\rho}(x) \; = \; \frac{\Ga(\a +\rho x)}{\Ga(\rho x)} & = & \frac{\a}{\rho\Ga(1-\a)}\int_0^\infty (1- e^{-x u}) \frac{e^{-u\rho^{-1}}}{(1-e^{-u\rho^{-1}})^{\a +1}}\, du,
\end{eqnarray*}
we obtain as above
\begin{equation}
\label{Tzeta}
\frac{\Ga(\rho)}{\Ga(\rho+\a)}\; \T(1+\a\rho^{-1},\rho^{-1}, (1-\a)\rho^{-1})\;\elaw\;\int_0^\infty e^{-\zeta^{(\a,\rho)}_t}\, dt
\end{equation}
where $\lacc \zeta^{(\a,\rho)}_t\! , \; t\ge 0\racc$ is the subordinator having Laplace exponent
$$\esp\lcr e^{-\lbd \zeta^{(\a,\rho)}_t}\rcr\; = \; e^{-t \Psi_{\a,\rho}(\lbd)}.$$
This leads to
$$\Farl\; \elaw\; \F_\rho\;\times\; \lpa \lbd \int_0^\infty e^{-\zeta^{(\a,\rho)}_t}\, dt \rpa^{\frac{1}{\rho}}.$$
Let us again notice that the identification (\ref{Tzeta}) can be deduced from Corollary 5 in \cite{LS} is the case $q = -\rho -\a < -\a$ and ${\hat \rho} =0$ - see also Remark 10 therein. Observe also that the limiting cases $\a = 0$ and $\a =1$ are consistent, with respectively $\zeta^{(0,\rho)}_t = \sigma^{(0)}_t$ and $\zeta^{(1,\rho)}_t = \rho\, t.$ 

\medskip

(b) Since $\Psi_{1-\a,\rho}(x)\Phi_{\a,\rho}(x) = \rho x,$ we have
$$ n!\; =\;\rho^n  \, \times\, b_n(1-\a,\rho)\,\times\, a_n(\a,\rho)$$
for all $\a\in[0,1],\rho > 0$ and $n\ge 1.$ By moment determinacy, this implies the following factorization of the unit exponential law
\begin{eqnarray}
\label{LRho}
\L & \elaw & \lpa \int_0^\infty \rho \lpa 1+ \sigma_t^{(1-\a)}\rpa^{-\rho -1+\a} \, dt\rpa \times \lpa\int_0^\infty \lpa (1- \sga_t)_+\rpa^{\rho -\a} \, dt\rpa\nonumber\\
& \elaw & \lpa \int_0^\infty \lpa 1+ \rho^{-1}\sigma_t^{(1-\a)}\rpa^{-\rho -1+\a} \, dt\rpa \times \lpa\int_0^\infty \lpa (1- \rho^{-1}\sga_t)_+\rpa^{\rho -\a} \, dt\rpa
\end{eqnarray}
which is valid for all $\a\in[0,1]$ and $\rho > 0.$ For $\rho = \a,$ this factorization reads
$$\L\;\elaw\; \lpa \int_0^\infty \frac{\a\, dt}{1+ \sigma_t^{(1-\a)}}\rpa \times\; \Z_\a^{-\a}\; \elaw\; \L^\a\,\times\,\Z_\a^{-\a},$$
where the first identity follows from Remark \ref{rW1} (d) and the second one from (3.4) in \cite{LS}. The simple identity $\L\elaw \L^\a\times\Z_\a^{-\a}$ is well-known as Shanbhag-Sreehari's identity. It has been thoroughly discussed in Section 3 of \cite{BY2001} from the point of view of perpetuities of subordinators, and their associated remainders. Observe also that changing the variable and letting $\rho\to\infty$ in (\ref{LRho}) leads to
$$\L\;\elaw\; \lpa \int_0^\infty e^{-\sigma_t^{(1-\a)}} dt \rpa \times\lpa \int_0^\infty e^{-\sigma_t^{(\a)}} dt \rpa,$$
another classic identity obtained in \cite{CPY} - see Example E therein. Last, it is interesting to mention the following identity, which follows at once from (\ref{LRho}), Theorem \ref{Weib} and Theorem \ref{Fresh}:
$$\W_{\a,1,\rho}\,\times\,\F_{1-\a,\rho,\rho}^{-1}\;\elaw\;\lpa \frac{\L\times\L}{\L} \rpa^{\frac{1}{\rho}}.$$

(c) The above proof shows that the function
$$x\;\mapsto\;E_{\alpha,m,m-\frac{1}{\a}}(-x)$$
is CM for every $\a\in (0,1]$ and $m >0.$ This is a generalization of the fact that $ E_{\alpha,1,1-\frac{1}{\a}}(-x)\; =\; \Ga(\a)E_{\a,\a}(-x)$ is CM for every $\a\in (0,1],$ which is itself a direct consequence of the aforementioned Pollard theorem because $\a E_\a'(-x) = E_{\a,\a} (-x).$ The formula (\ref{KS2}) implies the Bernstein representation
$$E_{\alpha,m,m-\frac{1}{\a}}(-x)\; =\; \esp\lcr \exp -x \lacc \int_0^\infty \lpa 1+ \sga_t\rpa^{-\a (m+1)} \, dt\racc\rcr, \qquad x\ge 0.$$
For $m = 1,$ with the notation of Remark \ref{rW1} (d) we obtain
$$E_{\alpha,1,1-\frac{1}{\a}}(-x)\; =\; \Ga(\a)E_{\a,\a}(-x) \; =\; \Ga(1+\a)E'_{\a}(-x)\; =\; \Ga(1+\a)\esp \lcr T_\a\, e^{-x T_\a}\rcr\; =\; \esp \lcr e^{-x T_\a^{(1)}}\rcr$$
where $T^{(1)}_\a$ is the size-bias of order 1 of $T_\a.$ This implies the curious identity
$$T_\a^{(1)} \; \elaw\; \int_0^\infty \lpa 1+ \sga_t\rpa^{-2\a} \, dt.$$
In a different direction, is is worth recalling that for every $\beta > \a$ and $\a\in (0,1]$ the function
$$E_{\alpha,1,\frac{\beta-1}{\a}}(-x)\; =\; \Ga(\beta)E_{\a,\beta}(-x)\; =\;\esp \lcr e^{-x \,\B_{\a,\beta -\a}^\a \times\, T_\a^{(1)}}\rcr$$
is also CM, where the Bernstein representation involving the $\a-$power of a standard Beta distribution follows directly from Lemma 4.26 in \cite{GKMR} - see also the references therein. In Section \ref{FurthA} below, we will come back to this example together with further analytical properties of the Kilbas-Saigo functions $E_{\alpha,m,m-\frac{1}{\a}}.$}

\end{rem}

We end this section with a convergent series representation in the non-explicit case $\a\in(0,1)$ for the density $f^\F_{\a,\lbd,\rho}$ of $\Farl.$ This is a consequence of a term-by-term differentiation of the distribution function
$$\pb[\Farl \le x]\; =\; E_{\a,\frac{\rho}{\a},\frac{\rho-1}{\a}}(-\lbd x^{-\rho})$$
which was obtained during the proof of Theorem \ref{Fresh}.
 
\begin{coro}
\label{FreshD}
For every $\a\in (0,1),$ the density of $\Farl$ has the following convergent series representation on $(0,\infty)$
$$f^\F_{\a,\lbd,\rho}(x)\; =\; \lbd \, x^{-\rho -1}\,\sum_{n\ge 0}\lpa \prod_{j=0}^n \frac{\Ga(j\rho + \rho + 1)}{\Ga(j\rho +\rho +\a)}\rpa \frac{(-\lbd x^{-\rho})^n}{\rho^n\,n!}\cdot$$
\end{coro}

\subsection{Proof of Theorem \ref{Gumb}}

The argument is shorter than for Theorems \ref{Weib} and \ref{Fresh}. We first consider the case $\a\in (0,1).$ The uniqueness is a direct consequence of Theorem \ref{t1G}. We next compute, by Fubini's theorem, the survival function
$$\pb[\Gal > x]\; =\; \pb[\L > e^{\lbd x +\G_\a}] \; = \; \esp\lcr e^{-e^{\lbd x +\G_\a}}\rcr\; =\; \sum_{n\ge 0} \frac{(-1)^n e^{\lbd n x}}{n!}\, \esp\lcr e^{n\G_\a}\rcr$$
for every $x\in\rl.$ On the other hand, since
$$e^{\G_\a}\;\elaw\;\int_0^\infty e^{-\sga_t}\, dt,$$
we know from Proposition 3.3 in \cite{CPY} that $\esp\lcr e^{n\G_\a}\rcr = (n!)^{1-\a}$ for all $n\ge 0.$ 
This implies
\begin{equation}
\label{LREM}
\bF (x)\; =\; \pb[\Gal > x]\; =\;\sum_{n\ge 0} \frac{(-1)^n e^{\lbd n x}}{(n!)^\a}\; =\; \cL_\a(-e^{\lbd x}).
\end{equation} 
A direct integration based on the Gamma integral and Fubini's theorem shows finally that
$$\Dp F (x)\; =\; \lbd^\a e^{\lbd x} \bF(x)$$
for every $x\in\rl$ as required. The case $\a = 1$ was already discussed in the introduction with a unique solution $F(x) = \pb[\G \le \lbd x] = \pb[\G_{1,\lbd} \le x]$, whereas the unique solution in the case $\a = 0$ is obviously $F(x) = 1/(1+ e^{-\lbd x}),$ which is the distribution function of $\lbd^{-1} (\G -\G) = \G_{0,\lbd}.$

\qed

\begin{rem}
\label{rG1}
{\em (a) The above proof also shows that the unique distribution function solving the fractional differential equation
$$\Dm \bF (x)\; =\; \lbd^\a e^{-\lbd x} F(x), \qquad F(x) > 0\;\mbox{on $\rl,$}$$
is $F(x) = \cL_\a(-e^{-\lbd x}) = \pb[-\Gal \le x].$

\medskip

(b) It is easy to deduce from the representations of the fractional extreme distributions in terms of integrals of the stable subordinator the following convergences in law
$$\rho\lbd^{-1}\lpa\W_{\a, \rho^\a\!, \rho} - 1\rpa\;\claw\;\Gal\qquad\mbox{and}\qquad \rho\lbd^{-1}\lpa 1- \F_{\a,\rho^\a\! ,\rho}\rpa\;\claw\;\Gal$$
as $\rho\to\infty,$ for every $\a\in [0,1]$ and $\lbd > 0.$ Observe that the case $\a=\lbd =1$ amounts to the aforementioned convergences in law $\rho(\W_\rho -1)\claw\G$ and $\rho(1-\F_\rho)\claw\G.$}

\end{rem}

As above, we finish this paragraph with a convergent series representation in the non-explicit case $\a\in(0,1)$ for the density $f^\G_{\a,\lbd}$ of $\Gal,$ which is a consequence of a term-by-term differentiation of the survival function $\pb[\Gal > x] = \cL_{\a}(-e^{\lbd x}).$

\begin{coro}
\label{GumbD}
For every $\a\in (0,1],$ the density of $\Gal$ has the following convergent series representation on $\rl$
$$f^\G_{\a,\lbd}(x)\; =\; \sum_{n\ge 1} \frac{(-1)^{n-1} \lbd n\, e^{\lbd n x}}{(n!)^\a}\cdot$$
\end{coro}

\section{Further properties}

\label{FurthA}

\subsection{On the complete monotonicity of the Kilbas-Saigo function} In this paragraph, motivated by the previous examples arising from the fractional Weibull and Fr\'echet distributions, we wish to characterize the complete monotonicity of all functions $x\mapsto E_{\a,m,l}(-x)$ on $(0,\infty).$ We begin with the following result on generalized Pochhammer symbols, which is reminiscent of Proposition 5.1 and Theorem 6.2 in \cite{D2010} and has an independent interest.

\begin{lem}
\label{Gtype}
Let $a,b,c,d$ and $\delta$ be positive parameters. There exists a positive random variable $Z =\Z[a,c;b,d;\delta]$ such that 
\begin{equation}
\label{Gdelta}
\esp[Z^s]\; =\; \frac{[a;\delta]_s[c;\delta]_s}{[b;\delta]_s[d;\delta]_s}
\end{equation}
for every $s > 0,$ if and only if $b+d \le a+c$ and $\inf\{b,d\} \le \inf\{a, c\}.$ This random variable is absolutely continuous on $(0,\infty)$, except in the degenerate case $a=b=c=d.$ Its support is $[0,1]$ if $b+d = a+c$ and $[0,\infty)$ if $b+d < a+c.$ 
\end{lem}

\proof

We discard the degenerate case $a=b=c=d,$ which is obvious with $Z=1.$ By (\ref{Malm}) and some rearrangements - see also (2.15) in \cite{LS}, we first rewrite
$$\log\lpa \frac{[a;\delta]_s[c;\delta]_s}{[b;\delta]_s[d;\delta]_s}\rpa \; =\; \kappa\, s\; + \; \int_{-\infty}^0 (e^{sx} - 1 - sx) \lpa \frac{e^{- b\vert x\vert} + e^{- d\vert x\vert} - e^{- a\vert x\vert} - e^{- c\vert x\vert}}{\vert x\vert (1-e^{-\vert x\vert})(1-e^{-\delta \vert x\vert})}\rpa dx$$
for every $s > 0,$ where $\kappa$ is some real constant. By convexity, it is easy to see that if $b+d \le a+c$ and $\inf\{b,d\} \le \inf\{a, c\},$ then the function $z\mapsto z^b + z^d - z^a - z^c$ is positive on $(0,1).$ This implies that the function
$$x\;\mapsto\;  \frac{e^{- b\vert x\vert} + e^{- d\vert x\vert} - e^{- a\vert x\vert} - e^{- c\vert x\vert}}{\vert x\vert (1-e^{-\vert x\vert})(1-e^{-\delta \vert x\vert})}$$
is positive on $(-\infty, 0)$ and that it can be viewed as the density of some L\'evy measure on $(-\infty,0),$ since it integrates $1\wedge x^2.$ By the L\'evy-Khintchine formula, there exists a real infinitely divisible random variable $Y$ such that 
$$\esp[e^{sY}] \; =\;\frac{[a;\delta]_s[c;\delta]_s}{[b;\delta]_s[d;\delta]_s}$$  
for every $s > 0,$ and the positive random variable $Z = e^Y$ satisfies (\ref{Gdelta}). Since we have excluded the degenerate case, the L\'evy measure of $Y$ is clearly infinite and it follows from Theorem 27.7 in \cite{S1999} that $Y$ has a density and the same is true for $Z.$

Assuming first $b+d = a+c,$ a Taylor expansion at zero shows that the density of the L\'evy measure of $Y$ integrates $1\wedge \vert x\vert$ and we deduce from (\ref{Malm}) the simpler formula
$$\log\esp[e^{sY}]\; =\; \log\lpa \frac{[a;\delta]_s[c;\delta]_s}{[b;\delta]_s[d;\delta]_s}\rpa \; =\; -\int^{\infty}_0 (1 - e^{-sx}) \lpa \frac{e^{- bx} + e^{- dx} - e^{- ax} - e^{- cx}}{x(1-e^{-x})(1-e^{-\delta x})}\rpa dx.$$
By the L\'evy-Khintchine formula, this shows that the ID random variable $Y$ is negative. Moreover, its support is $(-\infty,0]$ since its L\'evy measure has full support and its drift coefficient is zero - see Theorem 24.10 (iii) in \cite{S1999}, so that the support of $Z$ is $[0,1].$ 

Assuming second $b+d < a+c,$ the same Taylor expansion as above shows that the density of the L\'evy measure of $Y$ does not integrate $1\wedge \vert x\vert$ and the real L\'evy process associated to $Y$ is hence of type C with the terminology of \cite{S1999} - see Definition 11.9 therein. By Theorem 24.10 (i) in \cite{S1999}, this implies that $Y$ has full support on $\rl,$ and so does $Z$ on $\rl^+\!.$ \\
  
It remains to prove the only if part of the Lemma. Assuming $a \le d$ and $b\le c$ without loss of generality, we first observe that if $a < b$ then the function
$$s\;\mapsto\; \frac{[a;\delta]_s[c;\delta]_s}{[b;\delta]_s[d;\delta]_s}$$
is real-analytic on $(-b,\infty)$ and vanishes at $s = - a > -b,$ an impossible property for the Mellin transform of a positive random variable. The necessity of $b+d\le a+c$ is slightly more subtle and hinges again upon infinite divisibility. First, setting $\varphi(z) = z^b + z^d - z^a - z^c$ and $z_* =\inf\{ z > 0, \; \varphi(z) < 0\},$ it is easy to see by convexity and a Taylor expansion at 1 that if $b+d > a+c,$ then $z_* < 1$ and $\varphi(z) < 0$ on $(z_*,1)$ with $\varphi(z)\sim (b+d-a-c)(z-1)$ as $z\to 1.$ Introducing next the ID random variable $V$ with Laplace exponent
$$\log\esp[e^{sV}]\; =\; - \kappa\, s\; + \; \int_{\log z_*}^0 (e^{sx} - 1 - sx) \lpa \frac{e^{- a\vert x\vert} + e^{- c\vert x\vert} - e^{- b\vert x\vert} - e^{- d\vert x\vert}}{\vert x\vert (1-e^{-\vert x\vert})(1-e^{-\delta \vert x\vert})}\rpa dx,$$
we obtain the decomposition
$$\log\lpa \frac{[a;\delta]_s[c;\delta]_s}{[b;\delta]_s[d;\delta]_s}\rpa\; +\; \log\esp[e^{sV}] \; =\;\int_{-\infty}^{\log z_*} (e^{sx} - 1 - sx) \lpa \frac{e^{- b\vert x\vert} + e^{- d\vert x\vert} - e^{- a\vert x\vert} - e^{- c\vert x\vert}}{\vert x\vert (1-e^{-\vert x\vert})(1-e^{-\delta \vert x\vert})}\rpa dx,$$
whose right-hand side is the Laplace exponent of some ID random variable $U$ having an atom because its L\'evy measure, whose support is bounded away from zero, is finite - see Theorem 27.4 in \cite{S1999}. On the other hand, the random variable $V$ has an absolutely continous and infinite L\'evy measure and hence it has also a density. If there existed $Z$ such that (\ref{Gdelta}) holds, then the independent decomposition $U\elaw V  + \log Z$ would imply by convolution that $U$ has a density as well. This contradiction finishes the proof of the Lemma.

\endproof

\begin{rem} {\em (a) By the Mellin inversion formula, the density of $\Z[a,c;b,d;\delta]$ is expressed as
$$f(x) \; =\; \frac{1}{2\ii\pi x} \int_{s_0 -\ii \infty}^{s_0 + \ii\infty} x^{-s} \lpa \frac{[a;\delta]_s[c;\delta]_s}{[b;\delta]_s[d;\delta]_s} \rpa \, ds$$
over $(0,\infty)$ for any $s_0 > -\inf\{b,d\}.$ From this expression, it is possible to prove that this density is real-analytic over the interior of the support. We omit details. Let us also mention by Remark 28.8 in \cite{S1999} that this density is positive over the interior of its support.

\medskip

(b) With the standard notation for Pochhammer symbols, the aforementioned Proposition 5.1 and Theorem 6.2 in \cite{D2010} show that
$$s\;\mapsto\;\frac{(a)_s(c)_s}{(b)_s(d)_s}$$
is the Mellin transform of a positive random variable if and only if $b+d \ge a+c$ and $\inf\{b,d\} \ge \inf\{a, c\}.$ This fact can be proved exactly as above, in writing
$$\log\lpa \frac{(a)_s(c)_s}{(b)_s(d)_s}\rpa \; =\; - \int^{\infty}_0 (1-e^{-sx}) \lpa \frac{e^{- ax} + e^{- c x} - e^{- b x} - e^{- d x}}{x(1-e^{-x})}\rpa dx.$$
This expression also shows that the underlying random variable has support $[0,1]$ and that it is absolutely continuous, save for $a+c = b+d$ where it has an atom at zero. We refer to \cite{D2010} for an exact expression of the density on $(0,1)$ in terms of the classical hypergeometric function.}

\end{rem}

We can now state the main result of this paragraph, which characterizes the CM property for $E_{\a,m,l}(-x)$ on $(0,\infty).$

\begin{prop}
\label{KSCM}
Let $\a, m > 0$ and $l > -1/\a.$ The Kilbas-Saigo function
$$x\;\mapsto\; E_{\a,m,l}(-x)$$
is {\em CM} on $(0,\infty)$ if and only if $\a\le 1$ and $l\ge m-1/\a.$ Its Bernstein representation is  
\begin{equation}
\label{BKS}
E_{\alpha,m,l}(-x)\; =\; \esp\lcr \exp -x \lacc \X_{\a,m,l}\;\times\int_0^\infty \lpa 1+ \sga_t\rpa^{-\a (m+1)} \, dt\racc\rcr
\end{equation}
with $\delta =1/\a m$ and $\X_{\a,m,l} = \Z[1+1/m, (\a l +1)\delta;1, 1/m + (\a l +1)\delta;\delta].$ 
\end{prop}

\proof
Assume first $l\ge m-1/\a$ and let 
$$\Y_{\a,m,l}\;=\; \X_{\a,m,l}\;\times\int_0^\infty \lpa 1+ \sga_t\rpa^{-\a (m+1)}.$$ 
By Lemma \ref{Gtype} and Proposition 2 in \cite{LS}, its Mellin transform is
\begin{eqnarray*}
\esp[(\Y_{\a,m,l})^s] & = &  \delta^{s}\, \frac{[1+\delta;\delta]_s[(\a l +1)\delta;\delta]_s}{[1;\delta]_s[1/m +(\a l +1)\delta;\delta]_s}\\
& = & \Ga(1+s)\,\times\frac{[(\a l +1)\delta;\delta]_s}{[1/m +(\a l+1)\delta;\delta]_s}
\end{eqnarray*}
where in the second equality we have used (\ref{Conca4}). By Fubini's theorem, the moment generating function of $\Y_{\a,m,l}$ reads
\begin{eqnarray*}
\esp[e^{z \Y_{\a,m,l}}] & = & \sum_{n\ge 0}\; \esp[(\Y_{\delta,m,\eta})^n]\,\frac{z^n}{n!}\\
& = & \sum_{n\ge 0} \lpa\frac{[(\a l +1)\delta;\delta]_n}{[1/m +(\a l +1)\delta;\delta]_n}\rpa z^n\\
& = &  \sum_{n\ge 0} \lpa\prod_{j=0}^{n-1}\frac{\Ga(\a(jm +l)+1)}{\Ga(\a(jm +l +1)+1)}\rpa z^n\; =\; E_{\a,m,l}(z)
\end{eqnarray*}
for every $z\ge 0,$ where in the third equality we have used (\ref{Conca1}) repeatedly. The latter identity is extended analytically to the whole complex plane and we get, in particular,
$$ E_{\a,m,l}(-x) \; =\; \esp[e^{-x \Y_{\a,m,l}}],\qquad x\ge 0.$$
This shows that $E_{\a,m,l}(-x)$ is CM with the required Bernstein representation.\\

We now prove the only if part. If $E_{\a,m,l}(-x)$ is CM, then we see by analytic continuation that $E_{\a,m,l}(z)$ is the moment generating function on $\CC$ of the underlying random variable $X,$ whose positive integer moments read 
$$\esp[X^n]\; =\; n!\,\times\lpa\prod_{j=0}^{n-1}\frac{\Ga(\a(jm +l)+1)}{\Ga(\a(jm +l +1)+1)}\rpa,\quad n\ge 0.$$ 
If $\a > 1,$ Stirling's formula implies $\esp[X^n]^{\frac{1}{n}}\to 0$ as $n\to\infty$ so that $X\equiv 0,$ a contradiction because $E_{\a,m,l}$ is not a constant. If $\a = 1$ and $l+1 < m,$ then
$$\esp[X^n]\; = \;\frac{n!}{(c)_n m^n}\; \sim\; \frac{n^{1-c}}{m^n}\quad \mbox{as $n\to\infty,$}$$ 
with $c = (l+1)/m \in (0,1).$ In particular, the Mellin transform $s\mapsto\esp[X^s]$ is analytic on $\{\Re(s)\ge 0\},$ bounded on $\{\Re(s) = 0\},$ and has at most exponential growth on $\{\Re(s) > 0\}$ because
$$\vert \esp[X^s]\vert\; \le\; \esp\lcr X^{\Re(s)}\rcr\; =\; \lpa \esp\lcr X^{[\Re(s)] +1}\rcr\rpa^{\frac{\Re(s)}{[\Re(s)] +1}}$$
by H\"older's inequality. On the other hand, the Stirling type formula (\ref{Stirling}) implies, after some simplifications, 
$$\delta^{-s} \frac{[1+\delta;\delta]_s[c;\delta]_s}{[1;\delta]_s[c+\delta;\delta]_s}\; = \; \delta^{-s} s^{1-c} (1+ o(1))\quad \mbox{as $\vert s\vert \to\infty$ with $\vert \arg s\vert < \pi$}$$
and this shows that the function on the left-hand side, which is analytic on $\{\Re(s)\ge 0\},$ has at most linear growth on $\{\Re(s) = 0\}$ and at most exponential growth on $\{\Re(s) > 0\}.$ Moreover, the above analysis clearly shows that
$$\esp[X^n] \; = \; \delta^{-n} \frac{[1+\delta;\delta]_n[c;\delta]_n}{[1;\delta]_n[c+\delta;\delta]_n}$$
for all $n\ge 0$ and by Carlson's theorem - see e.g. Section 5.81 in \cite{T1939}, we must have
$$\esp[X^s] \; = \; \delta^{-s} \frac{[1+\delta;\delta]_s[c;\delta]_s}{[1;\delta]_s[c+\delta;\delta]_s}$$
for every $s > 0,$ a contradiction since Lemma \ref{Gtype} shows that the right-hand side cannot be the Mellin transform of a positive random variable if $\eta < 0.$ The case $\a < 1$ and  $l+1/\a < m$ is analogous. It consists in identifying the bounded sequence
$$ \frac{1}{n!}\,\times\lpa\prod_{j=0}^{n-1}\frac{\Ga(\a(jm +l +1)+1)}{\Ga(\a(jm +l)+1)}\rpa$$
as the values at non-negative integer points of the function
$$\delta^{-s}\times \frac{[1;\delta]_s[1/m + (\a l +1)\delta;\delta]_s}{[1+\delta;\delta]_s[(\a l +1)\delta;\delta]_s}\; =\; \delta^{-s} e^{-(1-\a) s \ln (s) + \kappa s + O(1)}\quad\mbox{as $\vert s\vert \to \infty$ with $\vert \arg s\vert < \pi,$}$$
where the purposeless constant $\kappa$ can be evaluated from (\ref{Stirling}). On $\{\Re(s) \ge 0\},$ we see that this function has growth at most $e^{\pi (1-\a) \vert s\vert /2}$ and we can again apply Carlson's theorem. We omit details.

\endproof

\begin{rem}

\label{RCM}

{\em (a) When $m =1,$ applying (\ref{Conca1}) we see that the random variable $\X_{\a, 1,l}$ has Mellin transform
$$\esp[(\X_{\a,1,l})^s]\; =\; \frac{[2;\delta]_s[l+1/\a;\delta]_s}{[1;\delta]_s[1 + l + 1/\a;\delta]_s}\; =\; \frac{(\a)_{\a s}}{(\beta)_{\a s}}$$
with $\beta = 1+\a l\ge \a.$ This shows $\X_{\a, 1,l}\elaw \B_{\a,\beta -\a}^\a,$ and we recover the Bernstein representation of the CM function $\Ga(\beta) E_{\a,\beta}(-x)$ which was discussed in Remark \ref{rF1} (c). Notice also the very simple expression for the Mellin transform
$$\esp[(\Y_{\a,1,l})^s] \; =\; \frac{\Ga(1+\a l)\Ga(1+s)}{\Ga(1 +\a (l+s))}\cdot$$

(b) Another simplification occurs when $l+1/\a = k m$ for some integer $k\ge 1.$ One finds
$$\esp[(\X_{\a,m,km -1/\a})^s]\; =\; \frac{[k;\delta]_s[1 +1/m;\delta]_s}{[1;\delta]_s[k +1/m;\delta]_s}\; =\; \prod_{j=1}^{k-1} \lpa \frac{(\a j m)_{u}}{(\a(jm+1))_{u}}\rpa$$
for $u =\a m s \ge 0,$ which implies
$$\X_{\a, m, km - 1/\a}\; \elaw\; \lpa \B_{\a m, \a}\times\cdots\times \B_{\a m (k-1), \a}\rpa^{\a m}.$$
In general, the law of the absolutely continuous random variable $\X_{\a,m,l}$ valued in $[0,1]$ seems to have a complicated expression.

\medskip

(c) As seen during the proof, the random variable $\Y_{\a,m,l}$ defined by the Bersntein representation 
$$E_{\a, m, l}(-x) \; =\; \esp [e^{-x \Y_{\a, m, l}}]$$
has Mellin transform 
\begin{equation}
\label{MellY}
\esp[(\Y_{\a,m,l})^s]\; =\; \Ga(1+s)\, \times\frac{[(\a l +1)\delta;\delta]_s}{[1/m +(\a l+1)\delta;\delta]_s}
\end{equation}
with $\delta = 1/\a m,$ for every $s > -1.$ By Fubini's theorem, this implies
$$\int_0^\infty E_{\a, m, l}(-x)\, x^{s-1}\, dx \; = \;  \Ga(s)\,\esp [\Y_{\a,m,l}^{-s}]\; =\; \Ga(s)\Ga(1-s)\, \times\frac{[(\a l +1)\delta;\delta]_{-s}}{[1/m +(\a l+1)\delta;\delta]_{-s}}$$
for every $s\in (0,1).$ Notice that this formula, which seems unnoticed in the literature on the Kilbas-Saigo function, remains valid for $l \in(-1/\a, m-1/\a)$ by the analyticity of the map $l\mapsto E_{\a, m, l} (-x).$ For $m = 1,$ we recover from (\ref{Conca1}) the formula
$$\int_0^\infty E_{\a, \beta}(-x)\, x^{s-1}\, dx \; = \; \frac{1}{\Ga(\beta)} \int_0^\infty E_{\a,1, \frac{\beta-1}{\a}}(-x)\, x^{s-1}\, dx \; = \;\frac{\Ga(s)\Ga(1-s)}{\Ga(\beta - \a s)}$$ 
which is given in (4.10.3) of \cite{GKMR}, as a consequence of the Mellin-Barnes representation of $E_{\a,\beta}(z).$ Notice that there is no such Mellin-Barnes representation for $E_{\a,m,l}(z)$ in general.}

\end{rem}

\subsection{Asymptotic behaviour of the densities} 

\label{ABTD}

In this paragraph we study the behaviour of the density functions of the fractional Weibull and Fr\'echet distributions at both ends of their support. To this end, we evaluate the Mellin transforms of $\Warl$ and $\Farl.$ The case of the fractional Gumbel distribution requires different arguments and will be handled in Paragraph \ref{SCLR}.

\subsubsection{The Weibull case} As a consequence of Corollary \ref{WeibD}, we first obtain
$$f^\W_{\a,\lbd, \rho}(x) \; \sim\; \lpa\frac{\lbd\,\Ga(\rho +1-\a)}{\Ga(\rho)}\rpa x^{\rho -1}\quad\mbox{as $x\to 0.$}$$
The behaviour of the density at infinity is less immediate and we will need the exact expression of the Mellin transform of $\Warl,$ which has an interest in itself.

\begin{prop}
\label{MASW}
The Mellin transform of $\Warl$ is 
$$\esp\lcr\Warl^s\rcr\; =\; \lpa\frac{\rho^\a}{\lbd}\rpa^{\frac{s}{\rho}}\Ga(1+ s\rho^{-1})\,\times\,\frac{[\rho +(1-\a);\rho]_{-s}}{[\rho;\rho]_{-s}}$$
for every $s\in (-\rho,\rho).$ As a consequence, one has
$$f^\W_{\a,\lbd, \rho}(x)\; \sim\; \lpa\frac{\rho }{\lbd \Ga(1-\a)}\rpa x^{-\rho -1} \quad \mbox{as $x\to\infty.$}$$
\end{prop}

\proof
We start with a more concise expression of (\ref{MellY}) for $l = m-1,$ which is a direct consequence of (\ref{Conca4}):
$$\esp[(\Y_{\a,\frac{\rho}{\a},\frac{\rho}{\a} -1})^s]\; =\; \rho^{-s} \times\frac{[1+(1-\a)\rho^{-1};\rho^{-1}]_s}{[1;\rho^{-1}]_s}\cdot$$
By Theorem \ref{Weib}, we deduce
\begin{eqnarray*}
\esp\lcr\Warl^s\rcr & = &\esp\lcr\lpa \frac{\L}{\lbd \Y_{\a,\frac{\rho}{\a},\frac{\rho}{\a} -1}}\rpa^{\frac{s}{\rho}}\rcr\\
& = & \lpa\frac{\rho}{\lbd}\rpa^{\frac{s}{\rho}}\Ga(1+ s\rho^{-1})\,\times\, \frac{[1+(1-\a)\rho^{-1};\rho^{-1}]_{-s\rho^{-1}}}{[1;\rho^{-1}]_{-s\rho^{-1}}}\\
& = & \lpa\frac{\rho^\a}{\lbd}\rpa^{\frac{s}{\rho}}\Ga(1+ s\rho^{-1})\,\times\,\frac{[\rho +(1-\a);\rho]_{-s}}{[\rho;\rho]_{-s}}
\end{eqnarray*}
for every $s\in (-\rho,\rho)$ as required, where the third equality comes from (\ref{Conca3}). The asymptotic behaviour of the density at infinity is then a standard consequence of Mellin inversion. First, we observe from the above formula and (\ref{zeroes}) that the first positive pole of $s\mapsto \esp\lcr\Warl^s\rcr$ is simple and isolated in the complex plane at $s = \rho,$ with
\begin{eqnarray*}
\esp\lcr\W_{\a,\lbd,\rho}^s\rcr & \sim & \lpa\frac{\rho^\a}{\lbd}\rpa\,\times\,\frac{[\rho +(1-\a);\rho]_{-\rho}}{[\rho;\rho]_{-s}} \\
& \sim & \lpa\frac{\rho^{\rho +\a}}{\lbd}\rpa\,\times\,\frac{[\rho +(1-\a);\rho]_{-\rho}}{[2\rho;\rho]_{-\rho}}\,\times\,(\rho)_{-s}\\ 
& =& \lpa\frac{\rho}{\lbd \Ga(1-\a)}\rpa \times\, \Gamma (\rho-s)\;\sim\; -\lpa\frac{\rho}{\lbd \Ga(1-\a)}\rpa\times\, \frac{1}{s-\rho}
\end{eqnarray*}
as $s\to \rho,$ where the second asymptotics comes from (\ref{Conca4}) and the third equality from (\ref{Conca2}). Therefore, applying Theorem 4 (ii) in \cite{FGD1995} - beware the correction $(\log x)^k \to (\log x)^{k-1}$ to be made in the expansion of $f(x)$ therein, we obtain
$$f^\W_{\a,\lbd,\rho}(x)\; \sim\; \lpa\frac{\rho }{\lbd \Ga(1-\a)}\rpa x^{-\rho -1} \quad \mbox{as $x\to\infty$}$$
as required.
\endproof

\begin{rem}
\label{rW2}
{\em (a) Another proof of the asymptotic behaviour at infinity can be obtained from (\ref{JSW2}). By multiplicative convolution the latter implies, setting $f^\cG_{\a,\rho}$ for the density of $\cG(\rho +1-\a,1-\a),$
\begin{eqnarray*}
f^\W_{\a,\lbd,\rho}(x) & = & \lbd \, x^{\rho-1} \int_0^\infty f^\cG_{\a,\rho}(y)\,y^{-\rho}\, e^{-\frac{\lbd}{\rho}(\frac{x}{y})^\rho} dy \\
& = & \lpa\frac{\lbd}{\rho}\rpa^{\frac{1}{\rho}}  \int_0^\infty f^\cG_{\a,\rho}\lpa x(\rho \lbd^{-1} t)^{-\frac{1}{\rho}}\rpa\,t^{-\frac{1}{\rho}}\,e^{-t}\, dt\\
& \sim & \lpa\frac{\rho }{\lbd \Ga(1-\a)}\int_0^\infty t \,e^{-t}\, dt \rpa x^{-\rho -1}\; =\; \lpa\frac{\rho }{\lbd \Ga(1-\a)}\rpa x^{-\rho -1}
\end{eqnarray*}
as $x\to\infty,$ where for the asymptotics we have used the Proposition in \cite{JSW2018} and a direct integration. This argument does not make use of Mellin inversion and is overall simpler than the above, but it does not convey to the Fr\'echet case. 

\medskip

(b) The Mellin transform simplifies for $\a = 0$ and $\a = 1\! :$ we recover
$$\esp[\W_{0,\lbd,\rho}^s] \; = \; \lbd^{-\frac{s}{\rho}}\Ga(1+s\rho^{-1})\Ga(1-s\rho^{-1})\qquad\mbox{and}\qquad \esp[\W_{1,\lbd,\rho}^s] \; = \; \lpa\frac{\rho}{\lbd}\rpa^{\frac{s}{\rho}}\Ga(1+s\rho^{-1})$$
in accordance with the scaling property $\Warl \elaw \lbd^{-1/\rho}\W_{\a,1,\rho}$ and (\ref{WeibIDs}), where the first equality follows from (\ref{Conca1}) and (\ref{Rhorho}). The Mellin transform takes a simpler form in two other situations. 

\begin{itemize}

\item For $\rho=\a,$ we obtain from (\ref{MellY}), (\ref{Conca1}) and (\ref{Conca2})
$$\esp[(\Y_{\a,1,0})^s]\; =\; \frac{\Ga(1+s)}{\Ga(1+\a s)}\; =\; \esp[\Z_\a^{-\a s}],$$
in accordance with Remark \ref{rW1} (d). This yields
$$\W_{\a,\lbd,\a}\;\elaw\; \lpa \frac{\L}{\lbd \Y_{\a,1,0}}\rpa^{\frac{1}{\a}}\;\elaw\; \lbd^{-\frac{1}{\a}}\,\Z_\a\,\times\,\L^{\frac{1}{\a}},$$ 
an identity which was already discussed for $\lbd = 1$ in the introduction as the solution to (\ref{Barre}). The Mellin transform reads
$$\esp[\W_{\a,\lbd,\a}^s] \; = \; \lbd^{-\frac{s}{\a}}\,\frac{\Ga(1+s\a^{-1})\Ga(1-s\a^{-1})}{\Ga(1-s)}\cdot$$

\item For $\rho=1-\a,$  where we obtain from (\ref{Conca2}) 
$$\esp[\W_{1-\rho,\lbd,\rho}^s] \; = \; \lpa\frac{\rho}{\lbd}\rpa^{\frac{s}{\rho}}\,\frac{\Ga(1+s\rho^{-1})\Ga(\rho -s)}{\Ga(\rho)}\qquad\mbox{and}\qquad \W_{1-\rho,\lbd,\rho}\; \elaw\; \lpa\frac{\rho}{\lbd}\rpa^{\frac{1}{\rho}}\,\L^{\frac{1}{\rho}}\,\times\,\GG^{-1}_{\rho}$$
having denoted by $\GG_t,$ here and throughout, the standard Gamma random variable with parameter $t.$

\end{itemize}

\medskip

(c) Integrating the density and using $\pb[\Warl > x] = E_{\a,\frac{\rho}{\a}, \frac{\rho}{\a}-1}(-\lbd x^\rho),$ we obtain at once the following asymptotic behaviour at infinity for any $\a\in (0,1]$ and $m > 0:$ 
$$E_{\a,m, m-1}(-x)\;\sim\; \frac{1}{\Ga(1-\a)\, x}\qquad\mbox{as $x\to\infty.$}$$ 
This behaviour seems unnoticed in the literature on the Kilbas-Saigo function, and turns out to be the same as that of the classical Mittag-Leffler function $E_\a(-x)$ - see e.g. (3.4.15) in \cite{GKMR}. It is actually possible to get the behaviour of $E_{\a,m, l}(-x)$ at infinity for all $l>-1/\a$ with the help of the Mellin transform computed in Remark \ref{RCM} (c). Notice however that the first positive pole might not be simple, for example when $m > 1$ and $l = m-1 -1/\a.$ We shall not discuss this behaviour here, save for $l = m-1/\a$ in the framework of the Fr\'echet distribution - see Remark \ref{rF2} (c). 

\medskip

(d) The four examples discussed in (b) above have Mellin transform expressed in terms of the quotient of a finite number of Gamma functions, making it possible to use a Mellin-Barnes representation of the density in order to get a full asymptotic expansion at infinity. For example, using the standard notation of Definition C.1.1 in \cite{AAR1999}, one obtains from (1.8.28) in \cite{KST} the expansion
$$f^\W_{\a,\lbd,\a}(x)\; =\; \lbd\, x^{\a-1} \,E_{\a,\a}(-\lbd x^\a)\; \sim\; \sum_{n\ge 1} \frac{n\a\, x^{-1-n\a}}{\lbd^n\Ga(1-n\a)}$$
which is everywhere divergent. Unfortunately, the Mellin transform of $\Warl$ might have poles of variable order and there does not seem to exist any general formula for the full asymptotic expansion at infinity of the density of $\Warl$.}
\end{rem}

\subsubsection{The Fr\'echet case} As a consequence of Corollary \ref{FreshD}, we first have
$$f^\F_{\a,\lbd, \rho}(x) \; \sim\; \lpa\frac{\lbd\,\Ga(\rho +1)}{\Ga(\rho+\a)}\rpa x^{-\rho -1}\quad\mbox{as $x\to \infty.$}$$
The behaviour of the density at zero is less immediate and we will need, as above, the exact expression of the Mellin transform of $\Farl,$ whose strip of analyticity is larger than for $\Warl.$ 

\begin{prop}
\label{MASF}
The Mellin transform of $\Farl$ is 
$$\esp\lcr\Farl^s\rcr\; =\; \lpa\frac{\rho^\a}{\lbd}\rpa^{-\frac{s}{\rho}}\Ga(1- s\rho^{-1})\,\times\,\frac{[\rho +1;\rho]_{s}}{[\rho +\a;\rho]_{s}}$$
for every $s\in (-\rho-\a,\rho).$ As a consequence, one has
$$f^\F_{\a,\lbd, \rho}(x)\; \sim\; \lpa\frac{\rho^{\frac{\a^2}{\rho}}(\rho +\a) }{\lbd^{1 +\frac{\a}{\rho}}}\,\Ga(1+\a)\, G(1-\a; \rho)\, G(1+\a; \rho)\rpa x^{\rho +\a -1} \quad \mbox{as $x\to 0.$}$$
\end{prop}

\proof The evaluation of the Mellin transform is done as for the fractional Weibull distribution, starting from the more concise expression
$$\esp[(\Y_{\a,\frac{\rho}{\a},\frac{\rho -1}{\a}})^s]\; =\; \rho^{-s} \times\frac{[1+\rho^{-1};\rho^{-1}]_s}{[1 +\a\rho^{-1};\rho^{-1}]_s}$$
and writing
$$\esp\lcr\Farl^s\rcr \; = \; \esp\lcr\lpa \frac{\L}{\lbd \Y_{\a,\frac{\rho}{\a},\frac{\rho-1}{\a}}}\rpa^{-\frac{s}{\rho}}\rcr\; =\;\lpa\frac{\rho^\a}{\lbd}\rpa^{-\frac{s}{\rho}}\Ga(1- s\rho^{-1})\,\times\,\frac{[\rho +1;\rho]_{s}}{[\rho +\a;\rho]_{s}}\cdot$$
The asymptotic behaviour of $f^\F_{\a,\lbd, \rho}(x)$ at zero follows then as for that of $f^\W_{\a,\lbd, \rho}(x)$ at infinity, in considering the residue at the first negative pole $s = -(\rho+\a)$ which is simple and isolated in the complex plane, applying Theorem 4 (i) in \cite{FGD1995} - with the same correction as above, and making various simplifications. We omit details. 

\endproof

\begin{rem}
\label{rF2}{\em (a) Comparing the Mellin transforms, Propositions \ref{MASW} and \ref{MASF} imply the interesting factorization
\begin{equation}
\label{FacWF}
\Warl^{-1}\;\elaw\; \Farl\;\times\;\Z(\rho +1-\a,\rho+\a;\rho,\rho+1;\rho).
\end{equation}
In general, it follows from Proposition \ref{KSCM} that for every $\a, m, \lbd > 0$ and $l > m- 1/\a,$ there exists a positive random variable having distribution function $E_{\a,m, l}(-\lbd x^{-\a m}),$ and which is given by (\ref{MellY}), (\ref{BKS}) and Theorem \ref{Fresh} as the independent product
$$\F_{\a,\lbd, \a m}\;\times\;(\X_{\a,m,l})^{\frac{1}{\a m}} \;\elaw\; \F_{\a,\lbd, \a m}\;\times\;\Z(\a l +1,\a (m+1);\a m,\a(l+1) +1;\a m),$$
where the identity in law follows from (\ref{Conca3}). In this respect, the fractional Fr\'echet distributions can be viewed as the ``ground state'' distributions associated to the Kilbas-Saigo functions $E_{\a,m, l},$ in the limiting case $l = m-1/\a.$ 

\medskip

(b) As above, the Mellin transform simplifies for $\a = 0,1\! :$ in accordance with (\ref{FreshIDs}), we get
$$\esp[\F_{0,\lbd,\rho}^s] \; = \; \lbd^{\frac{s}{\rho}}\,\Ga(1+s\rho^{-1})\Ga(1-s\rho^{-1})\qquad\mbox{and}\qquad \esp[\F_{1,\lbd,\rho}^s] \; = \; \lpa\frac{\lbd}{\rho}\rpa^{\frac{s}{\rho}}\Ga(1-s\rho^{-1}).$$
The Mellin transform also takes a simpler form in the same other situations as above.

\begin{itemize}

\item For $\rho=\a,$ with
$$\esp[\F_{\a,\lbd,\a}^s] \; = \; \lbd^{\frac{s}{\a}}\,\frac{\Ga(\a)\Ga(1+s\a^{-1})\Ga(1-s\a^{-1})}{\Ga(\a+s)}\cdot$$
This yields the identity $\F_{\a,\lbd,\a}\,\elaw\, \lbd^{\frac{1}{\a}}\,(\Z^{-1}_\a)^{(\a)}\,\times\,\L^{-\frac{1}{\a}},$ which was discussed in the introduction for $\lbd =1$ as the solution to (\ref{Barres}). This is also in accordance with Remark \ref{rF1} (c), since 
$$(T_\a^{(1)})^{\frac{1}{\a}}\;\elaw\; ((\Z^{-\a}_\a)^{(1)})^{\frac{1}{\a}}\;\elaw\;(\Z^{-1}_\a)^{(\a)}.$$ 
Notice that the constant appearing in the asymptotic behaviour of the density at zero is also simpler: one finds
\begin{equation}
\label{Fala}
f^\F_{\a,\lbd, \a}(x)\; \sim\; \lpa\frac{2\a\,\Ga(1+\a)}{\lbd^2\,\Ga(1-\a)}\rpa x^{2\a -1} \quad \mbox{as $x\to 0.$}
\end{equation}

\item For $\rho=1-\a,$ with
$$\esp[\F_{1-\rho,\lbd,\rho}^s] \; = \; \lpa\frac{\lbd}{\rho}\rpa^{\frac{s}{\rho}}\,\Ga(1-s \rho^{-1})\Ga(1+s)\quad\mbox{and}\quad\F_{1-\rho,\lbd,\rho}\; \elaw\; \lpa\frac{\lbd}{\rho}\rpa^{\frac{1}{\rho}}\,\L^{-\frac{1}{\rho}}\,\times\,\L.$$
Here, the density converges at zero to a simple constant: one finds
$$f^\F_{1-\rho,\lbd, \rho}(x)\; \to\; \lpa\frac{\rho}{\lbd}\rpa^{\frac{1}{\rho}} \Ga(1+\rho^{-1})\qquad \mbox{as $x\to 0.$}$$
\end{itemize}

\medskip

(c) Integrating the density and using $\pb[\Farl \le x] = E_{\a,\frac{\rho}{\a}, \frac{\rho -1}{\a}}(-\lbd x^{-\rho}),$ we obtain the following asymptotic behaviour at infinity for any $\a\in (0,1]$ and $m > 0:$ 
$$E_{\a,m, m-\frac{1}{\a}}(-x)\;\sim\; (\a m)^{\frac{\a}{m}}\Ga(1+\a)\, G(1-\a; \a m)\, G(1+\a; \a m) \, x^{-1-\frac{1}{m}}\qquad\mbox{as $x\to\infty.$}$$ 
For $m=1,$ this behaviour matches the first term in the full asymptotic expansion
$$E_{\a,1, 1-\frac{1}{\a}}(-x)\; =\; \Ga(\a) E_{\a,\a}(-x)\;\sim\; \Ga(\a)\sum_{n\ge 1} \frac{(-1)^n}{\Ga(-\a n) x^{n+1}}\cdot$$
As for $E_{\a,m, m-1}(-x),$ such a full asymptotic expansion seems difficult to obtain for all values of $m.$}
\end{rem}

\subsection{Optimal bounds for the distribution functions}

In Theorem 4 of \cite{TS2014}, the following uniform hyperbolic bounds are obtained for the classical Mittag-Leffler function:
\begin{equation}
\label{MLB}
\frac{1}{1 +\Ga(1-\a) x}\; \le\; E_\a(-x)\; \le\; \frac{1}{1 +\frac{1}{\Ga(1+\a)}\, x}
\end{equation}
for every $\a\in[0,1]$ and $x\ge 0.$ The constants in these inequalities are optimal because of the asymptotic behaviours
$$E_\a(-x)\; \sim\;\frac{1}{\Ga(1-\a) x}\quad\mbox{as $x\to\infty$}\qquad\mbox{and}\qquad 1- E_\a(-x)\; \sim\; \frac{x}{\Ga(1+\a)} \quad\mbox{as $x\to 0.$}$$
In this paragraph, we shall obtain analogous bounds for the Kilbas-Saigo functions $E_{\a,m,m-1}(-x)$ and $E_{\a,m,m-\frac{1}{\a}}(-x),$ which are associated to the fractional Weibull and Fr\'echet distributions. We begin with the following monotonicity properties, of independent interest.

\begin{prop}
\label{Mono}
Fix $\a\in (0,1]$ and $x\in\rl.$ The functions 
$$m\; \mapsto\; E_{\a,m,m-1}(x)\qquad\mbox{and}\qquad m\; \mapsto\; E_{\a,m,m-\frac{1}{\a}}(x)$$
are decreasing on $(0,\infty)$ if $x > 0$ and increasing on $(0,\infty)$ if $x < 0.$  
\end{prop}

\proof

 This is a direct consequence of (\ref{KS1}), (\ref{KS2}), and the fact that $\sga_t > 0$ for every $t > 0.$

\endproof

\begin{rem}
{\em It would be interesting to know if the same property holds for $m \mapsto E_{\a,m,m - l}(x)$ and any $l \le 1/\a.$ However, only the cases $l = 1$ and $l =1/\a$ seem to involve the $\a-$stable subordinator in a direct way.}
\end{rem} 

As in \cite{TS2014}, our analysis to obtain the uniform bounds will use some notions of stochastic ordering. Recall that if $X,Y$ are real random variables such that $\esp[\varphi(X)] \le \esp[\varphi(Y)]$ for every $\varphi :\rl \to\rl$ convex, then $Y$ is said to dominate $X$ for the convex order, a property which we denote by $X\prec_{cx} Y.$ The following result on convex orderings for infinite Beta products has an independent interest.

\begin{lem}
\label{Tcx}
For every $a,b,c > 0$ and $d \ge c,$ one has
$$\T(a,b,c)\;\prec_{cx}\; \T(a,b,d).$$ 
\end{lem}

\proof

By the definition of $\T(a,b,c)$ and the stability of the convex order by mixtures - see Corollary 3.A.22 in \cite{SSh}, it is enough to show
$$(a+b) \B_{a,b}\;\prec_{cx}\;(a+c) \B_{a,c}$$ 
for every $a,b > 0$ and $c \ge b.$ Using again Corollary 3.A.22 in \cite{SSh} and the standard identity $\B_{a,c}\,\elaw\,\B_{a,b}\times\B_{a+b,c-b},$ we are reduced to
$$\lpa\frac{a+b}{a+c}\rpa\; =\; \esp[\B_{a+b,c-b}]\;\prec_{cx}\;\B_{a+b,c-b}$$  
which is a direct consequence of Jensen's inequality.

\endproof

The following result is a generalization of (\ref{MLB}), which deals with the case $m=1$ only, to all Kilbas-Saigo functions $E_{\a,m,m-1}.$ The argument is much simpler than in the original proof of (\ref{MLB}).

\begin{prop}
\label{KSO1}
For every $\a\in[0,1], m >0$ and $x\ge 0,$ one has
$$\frac{1}{1 +\Ga(1-\a) x}\; \le\; E_{\a,m,m-1}(-x)\; \le\; \frac{1}{1 +\frac{\Ga(1+\a(m-1))}{\Ga(1+\a m)}\, x}\cdot$$
\end{prop}

\proof

The first inequality is a direct consequence of Proposition \ref{Mono}, which implies
\begin{eqnarray*}
E_{\a,m,m-1}(-x) & \ge & \esp\lcr\exp \lacc -x\, \int_0^\infty \lpa 1- \sga_t\rpa_+^{-\a} \, dt\racc\rcr\\
& = & \esp\lcr e^{-x\,\Ga(1-\a)\,\L}\rcr\; = \; \frac{1}{1 +\Ga(1-\a) x}
\end{eqnarray*}
for $x\ge 0,$ where the first equality follows from Theorem (b) (ii) in \cite{LS}. For the second equality, we come back to the infinite product representation
\begin{eqnarray*}
\int_0^\infty \lpa 1- \sga_t\rpa_+^{\rho -\a} \, dt & \elaw & \frac{\Ga(\rho +1-\a)}{\Ga(\rho +1)}\; \T(1,\rho^{-1}, (1-\a)\rho^{-1})
\end{eqnarray*}
which follows from Theorem (b) (i) in \cite{LS}, as in the proof of Theorem \ref{Weib}. Lemma \ref{Tcx} implies then 
\begin{eqnarray*}
\int_0^\infty \lpa 1- \sga_t\rpa_+^{\rho -\a} \, dt & \prec_{cx} & \frac{\Ga(\rho +1-\a)}{\Ga(\rho +1)}\; \T(1,\rho^{-1}, \rho^{-1})\; \elaw\; \frac{\Ga(\rho +1-\a)}{\Ga(\rho +1)}\;\L
\end{eqnarray*}
where the identity in law follows from (2.7) in \cite{LS}. Using (\ref{KS1}) with $\rho = \a m$ and the convexity of $t\mapsto e^{-xt}$, we obtain the required
$$E_{\a,m,m-1}(-x)\; \le\; \frac{1}{1 +\frac{\Ga(1+\a(m-1))}{\Ga(1+\a m)}\, x}\cdot$$
\endproof

\begin{rem}
\label{rKS1}
{\em (a) As for the classical case $m=1$, these bounds are optimal because of the aforementioned asymptotic behaviours
$$E_{\a, m, m-1}(-x)\; \sim\;\frac{1}{\Ga(1-\a) x}\;\mbox{as $x\to\infty$}\quad\mbox{and}\quad 1\,-\, E_{\a, m, m-1}(-x)\; \sim\; \frac{\Ga(1+\a (m-1))}{\Ga(1+\a m)}\, x \;\mbox{as $x\to 0.$}$$

\medskip

(b) It is easy to check that the above proof also implies
$$E_{\a,m,m-1}(x)\; \le\; \frac{1}{(1 -\Ga(1-\a) x)_+}$$
for every $\a\in[0,1], m > 0$ and $x\ge 0.$ This seems unnoticed even in the classical case $m=1.$
}
\end{rem}

Our next result is a uniform hyperbolic upper bound for the Kilbas-Saigo function $ E_{\a, m, m-\frac{1}{\a}},$ with an optimal power exponent by Remark \ref{rF2} (c) and an optimal constant since 
$$1\,-\, E_{\a, m, m-\frac{1}{\a}}(-x)\; \sim\; \frac{\Ga(\a m)\, x}{\Ga(\a(m+1))}\quad \mbox{as $x\to 0.$}$$ 
\begin{prop}
\label{KSO2}
For every $\a\in(0,1], m >0$ and $x\ge 0,$ one has
$$E_{\a,m,m-\frac{1}{\a}}(-x)\; \le\; \frac{1}{\lpa 1 + \frac{\Ga(1+\a m)}{\Ga(1+\a (m+1))}\, x\rpa^{1+\frac{1}{m}}}\cdot$$ 
\end{prop}

\proof

The inequality is derived by convex ordering as in Proposition \ref{KSO1}: one has
\begin{eqnarray*}
\int_0^\infty \lpa 1+ \sga_t\rpa_+^{-\rho -\a} \, dt & \elaw & \frac{\Ga(\rho)}{\Ga(\rho +\a)}\; \T(1 +\a\rho^{-1},\rho^{-1}, (1-\a)\rho^{-1})\\
& \prec_{cx} & \frac{\Ga(\rho)}{\Ga(\rho +\a)}\; \T(1 +\a\rho^{-1},\rho^{-1}, \rho^{-1})\; \elaw\; \frac{\Ga(\rho +1)}{\Ga(\rho +1+\a)}\;\GG_{1+\frac{\a}{\rho}}
\end{eqnarray*}
where the first identity follows from Corollary 3 in \cite{LS} as in the proof of Theorem \ref{Weib}, the convex ordering from Lemma \ref{Tcx} and the second identity from (2.7) in \cite{LS}. Then, using (\ref{KS2}) with $\rho = \a m,$ we get the required inequality.
\endproof

As in Proposition \ref{KSO1} we believe that there is also a uniform lower bound, with a more complicated optimal constant which can be read off from the asymptotic behaviour of the density at zero obtained in Proposition \ref{MASF}:
\begin{equation}
\label{ConjLow}
E_{\a,m,m-\frac{1}{\a}}(-x) \;\ge\; \frac{1}{(1 +(\a m)^{-\frac{\a}{m+1}}(\Ga(1+\a)\, G(1-\a; \a m)\, G(1+\a; \a m))^{-\frac{m}{m+1}}\, x)^{1+\frac{1}{m}}}\cdot
\end{equation}
Unfortunately, the proof of this general inequality still eludes us. The monotonicity property observed in Proposition \ref{Mono} does not help here, giving only the trivial lower bound zero. The discrete factorizations which are used in \cite{TS2014} are also more difficult to handle in this context, because the Mellin transform underlying $E_{\a,m,m-\frac{1}{\a}}$ is expressed in terms of generalized Pochhammer symbols. In the case $m=1,$ we could however get a proof of (\ref{ConjLow}).

\begin{prop}
\label{KSO3}
For every $\a\in(0,1)$ and $x\ge 0,$ one has
$$E_{\a,1,1-\frac{1}{\a}}(-x)\; \ge\; \frac{1}{\lpa 1 +\sqrt{\frac{\Ga(1-\a)}{\Ga(1+\a)}}\, x\rpa^{\! 2}}\cdot$$
\end{prop}

\proof By Remark \ref{rF1} (c) we have $E_{\alpha,1,1-\frac{1}{\a}}(-x) = \esp \lcr e^{-x T_\a^{(1)}}\rcr$ and since
$$\esp \lcr e^{-x \,\GG_2}\rcr\; =\; \frac{1}{\lpa 1 + x\rpa^2}$$
for every $x \ge 0,$ as in the proof of Theorem 4 in \cite{TS2014} it is enough to show that
\begin{equation}
\label{StoO}
T_\a^{(1)}\,\prec_{st}\, \sqrt{\frac{\Ga(1-\a)}{\Ga(1+\a)}}\, \GG_2,
\end{equation}
where $\prec_{st}$ stands for the usual stochastic order between two real random variables. Recall that $X\prec_{st} Y$ means $\pb[X\ge x]\le\pb[Y\ge x]$ for every $x\in\rl.$ Since $T_{1/2}\elaw 2\sqrt{\GG_{1/2}},$ the case $\a =1/2$ is explicit and the stochastic ordering can be obtained directly. More precisely, the densities of both random variables in (\ref{StoO}) are respectively given by
$$\frac{x}{2}\; e^{-x^2/4}\qquad\mbox{and}\qquad \frac{x}{2}\; e^{-x/\sqrt{2}}$$
on $(0,\infty),$ where they cross only once at $x= 2\sqrt{2}.$ It is a well-known and easy result that this single intersection property yields (\ref{StoO}) - see Theorem 1.A.12 in \cite{SSh}.

The argument for the case $\a\neq 1/2$ is somehow analogous, but the details are more elaborate because the density of $T_\a^{(1)}$ is not explicit anymore. We proceed as in Theorem C of \cite{TS2014} and first consider the case where $\a$ is rational. Setting $\a = p/n$ with $n > p$ positive integers and $X_\a = T_\a^{(1)} = (Z_\a^{-\a})^{(1)}$ we have, on the one hand,
\begin{eqnarray*}
\esp[(X_{\a})^{n s}] & = & \frac{\esp[(T_\a)^{1+ns}]}{\esp[T_\a]} \\ 
& = & \frac{\Ga(2+ns) \Ga(1+pn^{-1})}{\Ga(1+pn^{-1} +ps)}\\
& = & \frac{n^{ns}}{p^{ps}}\,\times \,\esp\lcr \lpa\B_{\frac{2}{n}, \frac{1}{p} -\frac{1}{n}}\rpa^s\rcr\,\times\,\frac{\prod_{i=3}^{n+1} (in^{-1})_s}{\prod_{j=2}^{p} (jp^{-1} +n^{-1})_s}
\end{eqnarray*}
for every $s > -2n^{-1},$ where in the third equality we have used repeatedly the Legendre-Gauss multiplication formula for the Gamma function - see e.g. Theorem 1.5.2 in \cite{AAR1999}. The same formula implies, on the other hand,
\begin{eqnarray*}
\esp\lcr\lpa\sqrt{\frac{\Ga(1-\a)}{\Ga(1+\a)}}\; \GG_2\rpa^{\!\! ns}\, \rcr & = & \frac{n^{ns}\, \kappa_\a^s}{p^{ps}}\,\times \,\esp\lcr \lpa\GG_{\frac{2}{n}}\rpa^s\rcr\,\times\,\lpa \prod_{i=3}^{n+1} (in^{-1})_s\rpa\\
& = & \frac{n^{ns}}{p^{ps}}\,\times\,\esp\lcr \lpa\kappa_\a\,\times \,\GG_{\frac{2}{n}}\,\times\, \prod_{j=2}^{p}\; \GG_{\frac{j}{p} +\frac{1}{n}}\rpa^{\!s}\,\rcr\,\times\, \frac{\prod_{i=3}^{n+1} (in^{-1})_s}{\prod_{j=2}^{p} (jp^{-1} +n^{-1})_s}
\end{eqnarray*}
for every $s > -2n^{-1},$ with the notation
$$\kappa_\a\; =\; \lpa\prod_{i=1}^p\frac{\Ga(ip^{-1} - n^{-1})}{\Ga(ip^{-1} + n^{-1})}\rpa^{\frac{n}{2}}.$$ Since 
$$\frac{\prod_{i=3}^{n+1} (in^{-1})_s}{\prod_{j=2}^{p} (jp^{-1} +n^{-1})_s}\; =\;\esp\lcr \lpa\prod_{i=2}^{p}\; \B_{\frac{i+1}{n}, \frac{i}{p} - \frac{i}{n}}\,\times\! \prod_{j=p+1}^{n} \GG_{\frac{j+1}{n}}\rpa^{\!s}\;\rcr$$
for every $s > -3n^{-1},$ by factorization and Theorem 1.A.3(d) in \cite{TS2014} we are finally reduced to show
$$\B_{\frac{2}{n}, \frac{1}{p} -\frac{1}{n}} \prec_{st}  \lpa\prod_{i=1}^p\;\frac{\Ga(ip^{-1} - n^{-1})}{\Ga(ip^{-1} + n^{-1})}\rpa^{\frac{n}{2}} \!\times \,\GG_{\frac{2}{n}}\,\times\, \prod_{j=2}^{p}\; \GG_{\frac{j}{p} +\frac{1}{n}}$$
for every $n >p$ positive integers. The latter is equivalent to 
$$(\B_{\frac{2}{n}, \frac{1}{p} -\frac{1}{n}})^{\frac{2}{n}} \prec_{st}  \lpa\prod_{i=2}^p\;\frac{\Ga(ip^{-1} - n^{-1})}{\Ga(ip^{-1} + n^{-1})} \rpa\times \lpa\GG_{\frac{2}{n}}\,\times\, \prod_{j=2}^{p}\; \GG_{\frac{j}{p} +\frac{1}{n}}\rpa^{\frac{2}{n}}$$
and this is proved via the single intersection property exactly as for (5.1) in \cite{TS2014}: the random variable on the left-hand side has an increasing density on $(0,1),$ whereas the random variable on the right-hand side has a decreasing density on $(0,\infty),$ both densities having the same positive finite value at zero. We omit details. This completes the proof of (\ref{StoO}) when $\a$ is rational. The case when $\a$ is irrational follows then by a density argument.

\endproof

\begin{rem}
\label{rKS3}
{\em (a) It is easy to check from (\ref{Conca2}) and (\ref{Rhorho}) that
$$\frac{\Ga(1+\a)}{\Ga(1-\a)}\; =\;\a^\a\,\Ga(1+\a)\, G(1-\a; \a )\, G(1+\a; \a),$$
so that Proposition \ref{KSO3} leads to (\ref{ConjLow}) for $m=1,$ in accordance with the estimate (\ref{Fala}). In general, the absence of a tractable complement formula for the product $G(1-\a; \delta )\, G(1+\a; \delta)$ makes however the constant in (\ref{ConjLow}) more difficult to handle.

\medskip

(b) Combining Propositions \ref{KSO3} and \ref{KSO2} implies the following optimal bounds on the generalized Mittag-Leffler function $E_{\a,\a}(-x)$ for every $\a\in (0,1)$ and $x\ge 0,$ to be compared with (\ref{MLB}):
$$\frac{1}{\lpa 1 +\sqrt{\frac{\Ga(1-\a)}{\Ga(1+\a)}}\, x\rpa^2}\;\le\; \Ga(\a)\,E_{\a,\a}(-x)\; \le\; \frac{1}{\lpa 1 +\frac{\Ga(1+\a)}{\Ga(1+ 2\a)}\, x\rpa^2}\cdot$$
Notice that letting $\a\to 1$ leads to the trivial bound $0\le e^{-x}\le (2/(2+x))^2.$}
\end{rem}

Our last result in this paragraph gives analogous bounds for the generalized Mittag-Leffler functions $E_{\a,\beta}(-x)$ with $\a\neq\beta$ whenever they are completely monotone, that is for $\beta > \a$ - see Remark \ref{rF1} (c). Although there is no direct connection to fractional extreme laws, we include this result here because of its independent interest as a generalization of (\ref{MLB}).

\begin{prop}
\label{KSO4}
For every $\a\in(0,1], \beta > \a$ and $x\ge 0,$ one has the optimal bounds
$$\frac{1}{1 +\frac{\Ga(\beta-\a)}{\Ga(\beta)}\, x}\;\le\; \Ga(\beta)\,E_{\a,\beta}(-x)\; \le\; \frac{1}{1 +\frac{\Ga(\beta)}{\Ga(\beta+\a)}\, x}\cdot$$
\end{prop}

\proof By the last equality in Remark \ref{rF1} (c) we have
$$\Ga(\beta)\,E_{\a,\beta}(-x)\; =\; \esp\lcr e^{-x \, \Y_{\a,1,l}}\rcr$$
with $l = (\beta - 1)/\a > 1-1/\a$ and $\Y_{\a, 1, l} \elaw \B_{\a,\beta -\a}^\a \times\, T_\a^{(1)}.$ From Remark \ref{RCM} (a), one obtains
\begin{equation}
\label{MomXab}
\esp \lcr (\Y_{\a,1,l})^s\rcr\; =\; \frac{\Ga(1+s) \Ga(\beta)}{\Ga(\beta + \a s)}
\end{equation}
for every $s > -1,$ which implies the factorization $\L \elaw \Y_{\a,1,l}\times (\GG_\beta)^{\a}.$ Since, by Jensen's inequality,
$$\frac{\Ga(\beta +\a)}{\Ga(\beta)}\; =\; \esp\lcr (\GG_\beta)^\a\rcr \prec_{cx} (\GG_\beta)^{\a},$$
we deduce from Corollary 3.A.22 in \cite{SSh} the convex ordering
$$\Y_{\a,1,l}\prec_{cx} \frac{\Ga(\beta)}{\Ga(\beta +\a)}\;\L$$
which, as above, implies 
$$\Ga(\beta)\,E_{\a,\beta}(-x)\; \le\; \frac{1}{1 +\frac{\Ga(\beta)}{\Ga(\beta+\a)}\, x}$$
for every $x\ge 0.$ 

The argument for the other inequality is analogous to that of Proposition \ref{KSO3}. By density, we only need to consider the case $\a = p/n$ and $\beta = (p+q)/n$ with $p<n$ and $q$ positive integers. By (\ref{MomXab}) and the Legendre-Gauss multiplication formula, we obtain 
$$\esp \lcr (\Y_{\a,1,l})^{ns}\rcr\; =\;\frac{n^{ns}}{p^{ps}}\,\times \,\esp\lcr \lpa\B_{\frac{1}{n},\frac{q}{np}}\rpa^s\rcr\,\times\,\frac{\prod_{i=2}^{n} (in^{-1})_s}{\prod_{j=1}^{p-1} (jp^{-1} + (p+q)(np)^{-1})_s}$$
for every $s > -n^{-1}.$ On the other hand, one has
$$\esp\lcr \lpa \frac{\Ga(\beta-\a)}{\Ga(\beta)}\;\L\rpa^{\! ns}\,\rcr \; = \; \frac{n^{ns}}{p^{ps}}\;\esp\lcr \lpa\kappa_{\a,\beta}\,\times \,\GG_{\frac{1}{n}}\,\times\, \prod_{j=1}^{p-1}\; \GG_{\frac{j}{p} +\frac{p+q}{np}}\rpa^{\!\! s}\,\rcr\,\times\,\frac{\prod_{i=2}^{n} (in^{-1})_s}{\prod_{j=1}^{p-1} (jp^{-1} + (p+q)(np)^{-1})_s}$$
with
$$\kappa_{\a,\beta}\; =\; p^p \lpa\frac{\Ga(qn^{-1})}{\Ga((p+q)n^{-1})} \rpa^n.$$
Comparing these two formulas, we are reduced to show 
$$(\B_{\frac{1}{n},\frac{q}{np}})^{\frac{1}{n}}\,\prec_{st}\, p^{\frac{p}{n}} \lpa\frac{\Ga(qn^{-1})}{\Ga((p+q)n^{-1})} \rpa\times\lpa\GG_{\frac{1}{n}}\,\times\, \prod_{j=1}^{p-1}\; \GG_{\frac{j}{p} +\frac{p+q}{np}}\rpa^{\frac{1}{n}}$$
for every $p< n$ and $q$ positive integers. This is obtained in the same way as above via the single intersection property. We leave the details to the reader.

\endproof
 
\subsection{Some properties related to infinite divisibility} 

\label{SPID}

In this paragraph we derive some infinite divisibility properties for the fractional extreme distributions, in the spirit of the Corollary in \cite{JSW2018}. Recall that the law of a positive random variable $X$ is called a generalized Gamma convolution ($X\in\cG$ for short) if there exists a suitably integrable deterministic function $a : \rl^+\to\rl^+$ such that
$$X\; \elaw\; \int_0^\infty a(t)\, d\Gamma_t$$
where $\{\Ga_t, \, t\ge 0\}$ is the Gamma subordinator. Equivalently, one has $X\in\cG$ iff its log-Laplace exponent reads
$$-\log \esp[e^{-\lbd X}]\; =\; a\lbd\, +\, \int_0^\infty (1- e^{-\lbd x})\, k(x) \frac{dx}{x}$$
with $a \ge 0$ and $k(x)$ a CM function. This representation shows that a random variable $X\in\cG$ is also infinitely divisible ($X\in\cI$ for short). An important subclass of $\cG$ is that of hyperbolically completely monotone random variables, which we will denote by $\cH.$ By definition, one has $X\in\cH$ iff $X$ has a positive density $f_X$ on $(0,\infty)$ such that $f_X(uv) f_X(u v^{-1})$ is CM in the variable $v+v^{-1}$ for all $u >0.$ We refer to \cite{Bond} for a classic account on the classes $\cG$ and $\cH$, including the above facts and much more. See also \cite{JRY2008} for a more recent survey. In Chapter 7 of \cite{Bond}, the class $\cG$ is extended to distributions on the real line, under the denomination EGGC. More precisely, an infinitely divisible distribution on $\rl$ is called an EGGC if its L\'evy measure has a density on $\rl^*$ of the type $\vert x\vert^{-1} k(x)$ with $x\mapsto k(x)$ and $x\mapsto k(-x)$ being CM functions on $(0,\infty).$ In the following, we will say that such L\'evy measures belong to the Thorin class. We will also use the same notation $X\in\cG$ to denote EGGC distributions, since there shall be no ambiguity on the support of $X.$ 

Our analysis is based on the following two lemmas on infinite Beta products, which have an independent interest. The first one is a precision made on (2.5) in \cite{LS}, whereas the second one is an extension of the main argument for proof of the Corollary in \cite{JSW2018}. This extension was already discussed in Remark 1 therein but we give some detail for the sake of completeness, and for the independent interest of the logarithmic estimate (\ref{Tabc}).
 
\begin{lem}
\label{Pow}
For every $a,b,c > 0$ one has
$$\lpa \frac{\Ga(a+c)}{\Ga(a)\,b^c}\; \T(ab^{-1},b^{-1},cb^{-1})\rpa^{\!\frac{1}{b}}\,\elaw\,\;\frac{\Ga((a+c)b^{-1})}{\Ga(ab^{-1})}\; \T(a,b,c)$$
\end{lem}

\proof
By (2.5) in \cite{LS}, it is enough to show that the random variables on both sides have the same expectation. By Proposition 2 in \cite{LS}, this amounts to
$$\frac{G((a+c+1)b^{-1}; b^{-1})\, G(ab^{-1}; b^{-1})}{G((a+c)b^{-1}; b^{-1})\, G((a+1)b^{-1}; b^{-1})}\; =\; \frac{\Ga((a+c)b^{-1})\,b^{\frac{c}{b}}}{\Ga(ab^{-1})}$$
which is a consequence of (\ref{Conca2}). Alternatively, the identity in law can be obtained from (\ref{Conca3}).

\endproof

\begin{lem}
\label{HCMT}
For every $a,b,c > 0$ one has $\T(a,b,c)^{-1}\,\in\,\cG$ and 
$$\T(a,b,c)\,\in\,\cH\;\Leftrightarrow\;\T(a,b,c)\,\in\,\cI\;\Leftrightarrow\; c\,\ge\, b.$$
\end{lem}

\proof The fact that $\T(a,b,c)^{-1}\in\cG$ is derived as in Corollary 8 in \cite{LS}. For the if part of the equivalence, we first notice that the case $c=b$ is obvious since  by (2.7) in \cite{LS} one has $\T(a,b,b) \elaw a^{-1}\GG_a,$ whose density is HCM. If $c > b,$ we decompose
$$\T(a,b,c)\; \elaw\; \T(a,b,b)\,\times\,\T(a+b, b, c-b)\; \elaw\; a^{-1}\GG_a\,\times\, \T(a+b, b, c-b),$$
where the first identity in law follows from (2.4) in \cite{LS}, and we can conclude exactly as in the proof of the Corollary in \cite{JSW2018}, since $a+b >a.$ For the only if part of the equivalence, we need to show $c < b\Rightarrow \T(a,b,c)\not\in\cI.$ To do so, we first deduce from (\ref{Stirling}) and some simplifications the limit behaviour
$$\frac{1}{s} \lpa \esp\lcr \lpa \T(a,b,c)\rpa^{\frac{bs}{c}}\rcr\rpa^{\frac{1}{s}}\;\to \; \frac{1}{c \, {\rm e}}\,\lpa \frac{\Ga(ab^{-1})}{\Ga((a+c)b^{-1})}\rpa^{\frac{b}{c}}\qquad \mbox{as $s\to\infty.$}$$
Applying, as for the Corollary in \cite{JSW2018}, Lemma 3.2 in \cite{CSY1999}, implies
\begin{equation}
\label{Tabc}
\log \pb \lcr \T(a,b,c) > x\rcr\;\sim\; -\,c\, \lpa \frac{\Ga((a+c)b^{-1})}{\Ga(ab^{-1})}\, x\rpa^{\frac{b}{c}}\qquad \mbox{as $x\to\infty.$}
\end{equation}
This shows that $\T(a,b,c)\not\in\cI$ whenever $b >c,$ since its upper tail probabilities are then superexponentially small - see e.g. Theorem 26.8 in \cite{S1999}.

\endproof

\begin{rem}
\label{RTabc}
{\em For $a=b,$ it follows after some simplifications from the Theorem and the Proposition in \cite{JSW2018}, with $\a = c$ and $m = a+c$ therein, that the density of $\T(a,a,c)$ is equivalent to
$$\lpa\frac{(2\pi)^{\frac{a+c}{2} -1}\, a^{\frac{c}{a}}\,(\Ga(1+ca^{-1}))^{\frac{(a+c)(1+c)}{2c} -1}}{\sqrt{ac}\,G(a+c,a)}\rpa\, x^{\frac{(a+c)(1+c)}{2c} -2}\, e^{-c\,\lpa\Ga(1+ca^{-1})\, x\rpa^{\frac{a}{c}}}\qquad \mbox{as $x\to\infty.$}$$
This is an improvement at the natural scale of the logarithmic estimate (\ref{Tabc}). Notice that in the case $a=b=c,$ this behaviour matches the exact formula $\frac{a^{a} x^{a-1}}{\Ga(a)}\, e^{-ax}$ for the density of $\T(a,a,a)\elaw a^{-1}\GG_a.$ The exact behaviour of the density of $\T(a,b,c)$ at infinity for all $a,b,c > 0$ is an interesting open question.}
\end{rem}

We can now state the main result of this paragraph.

\begin{prop}
\label{HCMWFG}
For every $\lbd, \rho > 0$ and $\a\in[0,1]$ one has

\begin{itemize}

\item $\Warl\in\cG$ if $\rho\le 1$ and $\Warl\in\cH$ if $\rho\le 1-\a,$ 

\item $\Farl\in\cH$ if $\rho\le 1-\a,$

\item $\Gal\in\cG.$
\end{itemize}

\end{prop}

\proof

It follows from Theorem \ref{Weib}, Theorem \ref{Fresh}, Lemma \ref{Pow} and Formula (2.7) in \cite{LS} that 
$$\Warl\;\elaw\;\lpa\frac{\rho^\a}{\lbd}\rpa^{\frac{1}{\rho}} \frac{\Ga(1+\rho^{-1})}{\Ga(1+(1-\a)\rho^{-1})}\,\times\, \frac{\T(\rho,\rho,1)}{\T(\rho,\rho,1-\a)}$$and
$$\Farl\;\elaw\;\lpa\frac{\lbd}{\rho^\a}\rpa^{\frac{1}{\rho}} \frac{1}{\Ga(1+\a\rho^{-1})}\,\times\, \frac{\T(\rho+\a,\rho,1-\a)}{\T(\rho,\rho,1)}\cdot$$
If $\rho\le 1-\a,$ by the second statement in Lemma \ref{HCMT} all random variables involved on the right-hand side belong to $\cH$, and Theorem 5.1.1 in \cite{Bond} implies that $\Warl$ and $\Farl$ belong to $\cH$ as well. The fact that $\Warl\in\cG$ for $\rho\le 1$ follows from the first statement in Lemma \ref{HCMT} and Theorem 6.2.1 in \cite{Bond}. Last, a consequence of (\ref{Malm1}) is
\begin{eqnarray*}
\log \esp[e^{ s\Gal)}] & = & \log \Ga(1+ s\lbd^{-1})\; -\; (1-\a)\, \log \Ga(1+ s\lbd^{-1}) \\
& = & -\a\gamma\lbd^{-1} s\, +\, \int_\rl (e^{sx} - 1 - sx)\,(\Un_{\{x <0\}} \, +\, (1-\a)\Un_{\{ x > 0\}})\, \frac{dx}{\vert x\vert (e^{\lbd \vert x\vert} -1)}
\end{eqnarray*}
for every $s >-\lbd.$ By the L\'evy-Khintchine formula, this shows that $\Gal$ is infinitely divisible and that its L\'evy measure belongs to the Thorin class, in other words $\Gal\in\cG.$

\endproof

\begin{rem}
\label{RHCM}
{\em (a) The same proof shows that $\G_\a\in\cG,$ since
$$\log \esp[e^{ s\G_\a)}] \; = \; (1-\a)\, \log \Ga(1+ s) \; = \; (1-\a)\lpa -\gamma s\, +\, \int_{-\infty}^0 (e^{sx} - 1 - sx)\,\frac{dx}{\vert x\vert (e^{\vert x\vert} -1)}\rpa$$
for every $s > -1,$ so that $\G_\a\in\cI$ with a L\'evy measure in the Thorin class. On the other hand, Theorem 2.1 in \cite{BL2015} shows a superexponential behaviour for the density of $\L_\a$ at infinity, and this implies as above that $\L_\a\not\in\cI$ except for $\a = 0$ or $\a = 1.$ 

\medskip

(b) The above proof also shows that in general, $\Warl$ and $\Farl$ can be expressed as the independent quotient of two random variables in $\cG.$ Unfortunately, this does not allow one to infer further infinite divisibility properties. We believe however that $\Farl\in\cI$ for all $\rho > 0.$ See \cite{BS2013} for a panorama of results related to the infinite divisibility of the classical extreme distributions.}  
\end{rem}

\subsection{Some complements on the Le Roy function} 

\label{SCLR}

In this paragraph we derive some miscellaneous results on the Le Roy function
$$\cL_\a(x) \; =\; \sum_{n\ge 0} \frac{x^n}{(n!)^\a}\cdot$$
This function played a role in the proof of Theorem \ref{Gumb} and can be viewed as another generalization of the exponential function, which it is interesting to compare to the classical Mittag-Leffler function $E_\a(x).$ Here and throughout we discard the explicit cases $\cL_0(x) = E_0(x) = 1/(1-x)$ and $\cL_1(x) = E_1(x) = e^x.$\\

We begin with the asymptotic behaviour at infinity. Le Roy's original result - see \cite{L1899} p. 263 - reads 
$$\cL_\a (x)\; \sim\; \frac{(2\pi)^{\frac{1-\a}{2}}}{\sqrt{\a}}\, x^{\frac{1-\a}{2}}\, e^{\a x^{\frac{1}{\a}}}\qquad \mbox{as $x\to\infty,$}$$
and is obtained by a variation on Laplace's method. The latter method can be used to solve Exercise 8.8.4 in \cite{Olver}, which states
\begin{equation}
\label{Olv1}
\cL_\a (-x)\; = \; \frac{2 (2\pi)^{\frac{1-\a}{2}}}{\sqrt{\a}}\,  x^{\frac{1-\a}{2\a}}\, e^{\a \cos(\pi/\a) x^{\frac{1}{\a}}}\, \lpa \sin\lpa \pi (2\a)^{-1} + \a \sin\lpa\pi\a^{-1}\rpa x^{\frac{1}{\a}}\rpa\, +\, O(x^{-\frac{1}{\a}})\rpa
\end{equation}
for $\a\ge 2$ and 
\begin{equation}
\label{Olv2}
\cL_\a (-x)\; \sim \; \frac{1}{\a^\a\, \Ga(1-\a) \, x\, (\log x)^\a}
\end{equation}
for $\a \in (1,2),$ as $x\to \infty.$ The following estimate, which seems to have passed unnoticed in the literature, completes the picture.

\begin{prop}
\label{Olv3}
For every $\a\in (0,1),$ one has
$$\cL_\a (-x)\; \sim \; \frac{1}{\Ga(1-\a) \, x\, (\log x)^\a}\qquad \mbox{as $x\to\infty.$}$$
\end{prop}

\proof

By (\ref{LREM}), we have
$$\cL_\a (-x)\; =\; \pb\lcr \L > x \L_\a\rcr\; =\; \int_0^\infty e^{-xt} \, f_\a(t)\, dt$$
with $\L_\a = e^{\G_\a}$ having density $f_\a$ on $(0,\infty).$ On the one hand, recalling
$$\esp \lcr \L_\a^s\rcr\; =\; \Ga(1+s)^{1-\a}$$
for every $s > -1,$ we have $f_\a = e_{1-\a}$ with the notation of \cite{BL2015} and we can apply Theorem 2.4 therein to obtain
\begin{equation}
\label{Asfa}
f_\a(x)\;\sim\; \frac{1}{\Ga(1-\a)\, (-\log x)^\a}\qquad \mbox{as $x\to 0.$}
\end{equation}
Plugging this estimate into the above expression for $\cL_\a (-x),$ we conclude the proof by a direct integration.
\endproof

\begin{rem}
\label{Rolv}
{\em (a) The estimate (\ref{Asfa}) also gives the asymptotic behaviour of the density of $\Gal$ at the right end of the support. Indeed, by multiplicative convolution the density of $e^{\lbd \Gal}$ on $(0,\infty)$ writes
$$\int_0^\infty e^{-xy} \; y f_\a(y)\, dy \; \sim\; \frac{1}{\Ga(1-\a)\, x^2\, (\log x)^\a}\qquad\mbox{as $x\to\infty,$}$$
where the estimate follows from (\ref{Asfa}) as in the above proof. A change of variable implies then
$$f^\G_{\a,\lbd}(x)\;\sim\;\lpa\frac{\lbd^{1-\a}}{\Ga(1-\a)}\rpa\, x^{-\a}\, e^{-\lbd x}\qquad\mbox{as $x\to\infty.$}$$
The asymptotic behaviour of the density at the left end can be obtained as in Paragraph \ref{ABTD} via the moment generating function
$$\esp\lcr e^{s\Gal}\rcr\; =\; \Ga(1+ s\lbd^{-1})\,\Ga(1-s\lbd^{-1})^{1-\a}, \qquad \vert s\vert <\lbd.$$
Reasoning as in Proposition \ref{MASW} via the converse mapping theorem leads to
$$f^\G_{\a,\lbd}(-x)\;\sim\;\lbd\, e^{-\lbd x}\qquad\mbox{as $x\to\infty,$}$$
in accordance with the first term in the expansion given in Corollary \ref{GumbD}. Observe that this converse mapping argument does not work directly for estimating $f^\G_{\a,\lbd}(x)$ at the right end, because of the fractional singularity in the moment generating function. 

\medskip
 
(b) In the case $\a =2,$ one has $\cL_2(x) = I_0(2\sqrt{x})$ and $\cL_2(-x) = J_0(2\sqrt{x})$ for all $x\ge 0,$ where $I_0$ and $J_0$ are the classical, modified or not, Bessel functions with index 0. In particular, a full asymptotic expansion for $\cL_2$ at both ends of the support is available, to be deduced e.g. from (4.8.5) and (4.12.7) in \cite{AAR1999}. These expansions also exist when $\a$ is an integer since $\cL_\a$ is then a generalized Wright function - see Chapter F.2.3 in \cite{GKMR} and the original articles by Wright quoted therein. The case when $\a$ is not an integer is open, and might be technical in the absence of a true Mellin-Barnes representation.}
\end{rem}

Our next result characterizes the connection between the entire function $\cL_\a(z)$ and random variables. Recall that a function $f : \CC\to\CC$ which is holomorphic in a neighbourhood $\Omega$ of the origin is a moment generating function (MGF) if there exists a real random variable $X$ such that 
$$f(z) \; =\; \esp\lcr e^{z X}\rcr,\qquad z\in\Omega.$$
In particular, $\cL_0$ is the MGF of the exponential law $\L$ and $\cL_1$ is that of the constant variable $\Un.$ 

\begin{prop}
\label{LRCM}
The function $\cL_\a(z)$ is the {\em MGF} of a real random variable if and only if $\a\le 1.$ In this case, one has
$$\cL_\a(z)\; =\; \esp\lcr e^{z \L_\a}\rcr, \qquad z\in\CC.$$
\end{prop}

\proof

The if part is a consequence of (\ref{LREM}) as in the proof of Proposition \ref{Olv3}. For the only if part, the estimates (\ref{Olv1}) and (\ref{Olv2}) show that $\cL_\a(z)$ takes negative values on $\rl^-$ when $\a  > 1,$ so that it cannot be the moment generating function of a real random variable.

\endproof

Observe that since $\L_\a$ is non-negative, the above result also shows $\cL_\a(-x)$ is CM on $(0,\infty)$ if and only if $\a\le 1,$ echoing Pollard's classical result for the Mittag-Leffler $E_\a(-x)$ - see Proposition 3.23 in \cite{GKMR}. One can ask if there are further complete monotonicity properties for $\cL_\a,$ as in \cite{TS2015} for $E_\a.$ In a different direction, the following result gives a monotonicity property which is akin to Proposition \ref{Mono}.

\begin{prop}
\label{LRMono}
The mapping $\a\mapsto\cL_\a (x)$ decreases on $[0,1]$ for every $x\in\rl.$ 
\end{prop}
 
\proof The fact that $\a\mapsto\cL_\a (x)$ decreases on $\rl^+$ is obvious for $x\ge 0,$ by the definition of $\cL_\a.$ To show the property on $[0,1]$ for $x <0,$ we will use a convex ordering argument. More precisely, the Malmsten formula (\ref{Malm1}) and the L\'evy-Khintchine formula show that for every $t \in [0,1],$ the random variable $\G_{1-t} = \log \L_{1-t}$ is the marginal at time $t$ of a real L\'evy process, since $\esp[e^{\ii z \G_{1-t}}] =  \Ga(1+ \ii z)^t = e^{t \psi(z)}$ for every $z\in\rl,$ with
$$\psi(z)\; =\; -\gamma\ii z\, +\, \int^0_{-\infty} (e^{\ii zx} - 1 - \ii zx)\, \frac{dx}{\vert x\vert (e^{\vert x\vert} -1)}\cdot$$  
This is actually well-known - see again Example E in \cite{CPY}. By independence and stationarity of the increments of a L\'evy process, we deduce that there exists a multiplicative martingale $\{M_t, t\in [0,1]\}$ such that $M_t \elaw \L_{1-t}$ for every $t \in [0,1]$ and Jensen's inequality implies
$$\L_\beta\;\prec_{cx}\; \L_\a$$
for every $0\le\a\le\beta\le 1.$ Applying the definition of convex ordering to the function $\varphi(x) = e^{-x},$ we get $\cL_\beta(-x)\le\cL_\a(-x)$ for every $x > 0$ and $0\le\a\le\beta\le 1,$ as required. 
\endproof

\begin{rem}
\label{Picq}

{\em (a) In the terminology of \cite{HPRY}, the family $\{\L_{1-\a}, \; \a\in[0,1]\}$ is a peacock, whose associated multiplicative martingale is here completely explicit. We refer to \cite{HPRY} for numerous examples of explicit peacocks related to exponential functionals of L\'evy processes. Observe from Lemma \ref{Tcx} that the family $\{\T(a,b,t), \; t >0\}$ is also a peacock.

\medskip

(b) It is easily seen from Corollary 3 in \cite{LS} that $\T(1+\a q^{-1},q^{-1}, (1-\a) q^{-1})\claw \L_\a$ as $q\to\infty.$ On the other hand, we could not prove that the family $\{\T(1+(1-\a) q^{-1},q^{-1}, \a q^{-1}), \; \a\in[0,1]\}$ is a peacock for any fixed $q >0.$}
\end{rem}

Letting $\a\to 0,1$  in the above proposition leads to the bounds
$$e^{-x}\;\le\;\cL_\beta(-x)\;\le\;\cL_\a(-x)\;\le\; \frac{1}{1+x}$$
for every $x \ge 0$ and  $0\le\a\le\beta\le 1,$ to be compared with the less complete bounds (6.9) in \cite{TS2014} for $E_\a(-x).$ The hyperbolic upper bound is optimal as in Propositions \ref{KSO1} and \ref{KSO2}, since $1- \cL_\a(-x)\sim x$ as $x\to 0.$ The exponential lower bound is thinner than the order given in the estimate (\ref{Asfa}). On the other hand, it does not seem that stochastic ordering arguments can help to get a uniform estimate involving a logarithmic term. The following last proposition gives alternative bounds on $\cL_\a(-x)$ in terms of Kilbas-Saigo functions. It is a direct consequence of the Bernstein representations given in Remark \ref{rW1} (d), Remark \ref{rF1} (c) and Proposition \ref{LRCM}, and of the monotone character of the function $s\mapsto (1-xs^{-1})_+^s$ on $(0,\infty)$ for every $x\in\rl.$ We omit details. 

\begin{prop}
\label{LRKS}
For every $\a\in [0,1]$ and $x, m > 0,$ one has   
$$\lacc \begin{array}{l}
E_{\a, m,m-\frac{1}{\a}}(-(\a (m+1))^\a x)\;\le\;\cL_\a(-x)\;\le\; E_{\a, m+1,m}(-(\a m)^\a x),\\
E_{\a, m+1,m}((\a m)^\a x)\;\le\;\cL_\a(x)\;\le\; E_{\a, m,m-\frac{1}{\a}}((\a (m+1))^\a x).\end{array}\right.$$
Besides, $E_{\a, m,m-\frac{1}{\a}}((\a (m+1))^\a x)$ and $E_{\a, m+1,m}((\a m)^\a x)\to\cL_\a(x)$ as $m\to\infty,$ for every $x\in\rl.$ 
\end{prop}

\medskip

\section*{Appendix}

\subsection*{A.1. Fractional integrals and derivatives} In this paragraph, we fix the notation on the fractional operators which are used throughout the paper. This is an excerpt from the beginning of Chapter 2 in \cite{KST}. We will consider only three kinds of such operators which are the most familiar ones, and our fractional parameter $\a$ will always be supposed in $[0,1].$ There are certainly many other fractional operators with a larger family of fractional parameters, and we refer to the whole Chapter 2 in \cite{KST} for an account.

\subsubsection*{{\em A.1.1.} Progressive Liouville operators on the half-axis} For every $\a\in (0,1),$ the operator $f\mapsto \Iap(f)$ with 
$$\Iap (f)(x)\; =\; \frac{1}{\Ga (\a)}  \int_{0}^{x} (x-u)^{\alpha-1}f(u) \, du,\quad  x > 0,$$ 
is well-defined, taking possibly infinite values, on measurable functions $f : (0,\infty)\to \rl^+.$ It is easy to see that if $f$ is integrable at zero, then so is $\Iap(f)$ - see Lemma 2.1 in \cite{KST} for a more general result. The corresponding fractional derivative $f\mapsto \Dap f$ is such that
$$\Dap (f) (x)\; =\; \frac{\mathrm{d}}{\mathrm{d} x}\lpa {\rm I}^{1-\a}_{0+} (f)\rpa (x)$$
and is well-defined almost everywhere as soon as $f = \Iap(g)$ for some $g$ integrable at zero. Moreover, for such functions there is an inversion formula
 \begin{equation}
\label{Fr0+}
\Iap\lpa\Dap(f)\rpa\; =\; \Iap(g)\; =\; f,
\end{equation}
which is valid almost everywhere - see Lemma 2.5 in \cite{KST}. These operators are extended to the boundary cases $\a = 0$ with ${\rm I}^0_{0+} = {\rm D}^0_{0+} = {\rm Id}$ and $\a = 1$ with  ${\rm I}^1_{0+}$ and ${\rm D}^1_{0+}$ being respectively the usual running integral and derivative - see (2.1.7) in \cite{KST}. 
 
\subsubsection*{{\em A.1.2.} Regressive Liouville operators on the half-axis} 
 For every $\a\in (0,1),$ the operator $f\mapsto \Im(f)$ with 
$$\Im (f)(x)\; =\; \frac{1}{\Ga (\a)}  \int_{x}^\infty (u-x)^{\alpha-1}f(u) \, du,\quad  x > 0,$$ 
is well-defined on measurable functions $f : (0,\infty)\to \rl^+.$ It is easy to see that if $f$ is integrable at infinity, then so is $\Im(f)$ - see again Lemma 2.1 in \cite{KST}. The corresponding fractional derivative $f\mapsto \Dm f$ is such that
$$\Dm (f) (x)\; =\; -\frac{\mathrm{d}}{\mathrm{d} x}\lpa {\rm I}^{1-\a}_{-} (f)\rpa (x)$$
and is well-defined as soon as $f = \Im(g)$ for some $g$ integrable at infinity. Moreover, for such functions there is an inversion formula
 \begin{equation}
\label{Fr0-}
\Im\lpa\Dm(f)\rpa\; =\; \Im(g)\; =\; f,
\end{equation}
which is valid almost everywhere. These operators are extended to the boundary cases $\a = 0$ with ${\rm I}^0_{-} = {\rm D}^0_{-} = {\rm Id}$ and $\a = 1$ with  ${\rm I}^1_{-}$ and ${\rm D}^1_{-}$ being respectively the usual running integral and the opposite of the usual derivative - see again (2.1.7) in \cite{KST}. 

\subsubsection*{{\em A.1.3.} Progressive Liouville operators on the real axis}  For every $\a\in (0,1),$ the operator $f\mapsto \Ip(f)$ with 
$$\Ip (f)(x)\; =\; \frac{1}{\Ga (\a)}  \int_{-\infty}^x (x-u)^{\alpha-1} f(u) \, du,\quad  x \in\rl,$$ 
is well-defined on measurable functions $f : \rl \to \rl^+.$ Observe that $\Ip(f) (-x) = \Im (g) (x)$ with $g(x) = f(-x)$ so that we can transfer to $\Ip$ the properties on $\Im.$ In particular, if $f$ is integrable at $-\infty$, then so is $\Ip(f).$ The corresponding fractional derivative $f\mapsto \Dp f$ is such that
$$\Dp (f) (x)\; =\; -\frac{\mathrm{d}}{\mathrm{d} x}\lpa {\rm I}^{1-\a}_{+} (f)\rpa (x)$$
and is well-defined as soon as $f = \Ip(g)$ for some $g$ integrable at $-\infty$. Moreover, for such functions there is an inversion formula $\Ip\lpa\Dp(f)\rpa = \Ip(g) = f,$ which is valid almost everywhere. These operators are extended to the boundary cases $\a = 0$ with ${\rm I}^0_{+} = {\rm D}^0_{+} = {\rm Id}$ and $\a = 1$ with  ${\rm I}^1_{+}$ and ${\rm D}^1_{+}$ being respectively the usual running integral and derivative.

\subsection*{A.2. Barnes' double Gamma function} In this paragraph, mostly taken from \cite{BK1997} and \cite{K2011} to which we refer for further results, we gather some useful facts about Barnes' double Gamma function $G(z;\delta).$ For every $\delta > 0,$ this function is defined as the unique solution to the functional equation
\begin{equation}
\label{Conca1}
G(z+1;\delta)\; =\; \Gamma(z\delta^{-1}) G(z;\delta)
\end{equation}
with normalization $G(1;\delta) =1.$ This function is holomorphic on $\CC$ and admits the following Malmsten type representation
\begin{equation}
\label{Malm}
G(z;\delta)\; =\; \exp\int_0^\infty \lpa \frac{1-e^{-zx}}{(1-e^{-x})(1-e^{-\delta x})} \, -\, \frac{ze^{-\delta x}}{1-e^{-\delta x}}\, +\, (z-1)(\frac{z}{2\delta} -1)e^{-\delta x}\, - \, 1\rpa\frac{dx}{x}
\end{equation}
which is valid for $\Re(z) >0$ - see (5.1) in \cite{BK1997}. Putting (\ref{Conca1}) and (\ref{Malm}) together and making some simplifications, we recover the standard Malmsten formula for the Gamma function
\begin{equation}
\label{Malm1}
\Ga(1+z)\; =\; \exp\lacc -\gamma z\, +\, \int^0_{-\infty} (e^{zx} - 1 - zx)\, \frac{dx}{\vert x\vert (e^{\vert x\vert} -1)}\racc
\end{equation}
for every $z >-1,$ where $\gamma$ is Euler's constant. The following Stirling type asymptotic behaviour 
\begin{equation}
\label{Stirling}
\log G(z;\delta)\; -\; \frac{1}{2 \delta} \lpa z^2 \log z\, -\, (\frac{3}{2} + \log \delta)\, z^2 \, -\,(1+\delta)\, z\log z\rpa\; -\; A \, z\; -\; B \log z\; \to\; C
\end{equation}
is valid for $\vert z\vert\to\infty$ with $\vert \arg (z) \vert < \pi,$ for some real constants $A, B$ and $C$ which are given in (4.5) of \cite{BK1997}. There is a second concatenation formula
\begin{equation}
\label{Conca2}
G(z+\delta;\delta)\; =\; (2\pi)^{(\delta-1)/2} \delta^{1/2 -z}\Gamma(z) G(z;\delta)
\end{equation}
which is valid for all $z\in\CC,$ the right-hand side being understood as an anxiolytic extension when $z$ is a non-positive integer - see (4.6) in \cite{K2011} and the references therein. Observe that (\ref{Conca1}) and (\ref{Conca2}) lead readily to the closed formula
\begin{equation}
\label{Rhorho}
G(\delta;\delta)\; =\; G(1+\delta;\delta)\; =\; (2\pi)^{(\delta-1)/2} \delta^{-1/2}.
\end{equation}
In this paper we make an extensive use of the following Pochhammer type symbol
\begin{equation}
\label{Poch}
[a;\delta]_s\; =\; \frac{G(a +s;\delta)}{G(a;\delta)}
\end{equation}
which is well-defined for every $a,\delta > 0$ and $s > - a.$ The following formula
\begin{equation}
\label{Conca3}
[a\delta^{-1};\delta^{-1}]_{s\delta^{-1}} \; =\; (2\pi)^{s(1/\delta-1)/2}\, \delta^{s^2/2\delta - s(1+(1-2a)/\delta)/2}\,[a;\delta]_s
\end{equation}
can be deduced from (4.10) in \cite{K2011} - beware the different normalization for $G(1;\delta)$ therein which becomes irrelevant when considering the Pochhammer type symbol. Notice also that (\ref{Conca2}) yields
\begin{equation}
\label{Conca4}
\delta^s\,[a+\delta;\delta]_s\; =\; (a)_s \,
[a;\delta]_s
\end{equation}
with the standard notation
$$(a)_s \; =\; \frac{\Ga(a+s)}{\Ga(a)}$$
for the usual Pochhammer symbol. Finally, we observe from the double product representation of $G(z,\delta)$ - see e.g. (4.4) in \cite{K2011}, that for every $a, \delta > 0$ one has
\begin{equation}
\label{zeroes}
\inf\{s > 0, \;  [a;\delta]_{-s} = 0\}\; =\; a
\end{equation} 
and that this zero is simple and isolated on the complex plane.

\end{document}